\pdfoutput=1
\RequirePackage{ifpdf}
\ifpdf 
\documentclass[pdftex]{sigma}
\else
\documentclass{sigma}
\fi

\numberwithin{equation}{section}
\numberwithin{proposition}{section}
\numberwithin{lemma}{section}
\numberwithin{corollary}{section}
\numberwithin{definition}{section}

\newcommand{\coi}{C^\infty_0(\M)}
\newcommand{\A}{\mathcal{A}}
\newcommand{\Sc}{\mathcal{S}}
\newcommand{\HH}{\mathcal{H}}
\newcommand{\HHH}{\mathbb{H}}
\newcommand{\III}{\mathbb{I}}
\newcommand{\M}{\mathcal M}
\newcommand{\X}{\mathcal{X}}
\newcommand{\N}{\mathbb{N}}
\newcommand{\NN}{\mathcal{N}}
\newcommand{\R}{\mathbb{R}}
\newcommand{\C}{\mathbb{C}}
\newcommand{\de}{\mathrm{d}}
\newcommand{\abso}[1]{|#1|}
\newcommand{\norm}[1]{||#1||}
\newcommand{\dc}{d_D}

\newcommand{\pa}{{\mathcal{P}(\A)}}
\newcommand{\sa}{{\mathcal{S}(\A)}}
\newcommand{\tdx}{\tilde{\mathbf{x}}}

\begin{document}

\allowdisplaybreaks

\renewcommand{\thefootnote}{$\star$}

\renewcommand{\PaperNumber}{057}

\FirstPageHeading

\ShortArticleName{A View on Optimal Transport from Noncommutative Geometry}

\ArticleName{A View on Optimal Transport\\  from Noncommutative Geometry\footnote{This paper is a
contribution to the Special Issue ``Noncommutative Spaces and Fields''. The
full collection is available at
\href{http://www.emis.de/journals/SIGMA/noncommutative.html}{http://www.emis.de/journals/SIGMA/noncommutative.html}}}

\Author{Francesco D'ANDREA~$^\dag$ and Pierre MARTINETTI~$^\ddag$}

\AuthorNameForHeading{F.~D'Andrea and P.~Martinetti}

\Address{$^\dag$~Ecole de Math{\'e}matique, Univ.~Catholique de Louvain,\\
\hphantom{$^\dag$}~Chemin du Cyclotron 2, 1348 Louvain-La-Neuve, Belgium}
\EmailD{\href{mailto:francesco.dandrea@uclouvain.be}{francesco.dandrea@uclouvain.be}}

\Address{$^\ddag$~Institut f{\"u}r Theoretische Physik, Universit{\"a}t G{\"o}ttingen,\\
\hphantom{$^\ddag$}~Friedrich-Hund-Platz 1, 37077 G{\"o}ttingen, Germany}
\EmailD{\href{mailto:martinetti@theorie.physik.uni-goettingen.de}{martinetti@theorie.physik.uni-goettingen.de}}

\ArticleDates{Received April 14, 2010, in f\/inal form July 08, 2010;  Published online July 20, 2010}

\Abstract{We discuss the relation between the Wasserstein distance of order 1 between probability distributions on a metric space, arising in the study
of Monge--Kantorovich transport problem, and the spectral distance of noncommutative geometry. Starting from a~remark of Rief\/fel on compact manifolds, we f\/irst show that on any~-- i.e.~non-necessary compact~-- complete Riemannian spin manifolds, the two distances coincide. Then, on convex manifolds in the sense of Nash embedding, we provide some natural upper and lower bounds to the distance between any two probability distributions. Specializing
to the Euclidean space $\R^n$, we  explicitly compute the distance for a particular class of distributions genera\-li\-zing Gaussian wave packet. Finally we explore the analogy between the spectral and the Wasserstein distances in the noncommutative case, focusing on the standard model and the Moyal plane. In particular we point out that in the two-sheet space of the standard model, an optimal-transport interpretation of the metric requires a cost function that does not va\-nish on the diagonal. The latest is similar to the cost function occurring in the relativistic heat equation.}


\Classification{58B34; 82C70}

\renewcommand{\thefootnote}{\arabic{footnote}}
\setcounter{footnote}{0}

\section{Introduction}
The idea at the core of Noncommutative Geometry~\cite{Con94} is the observation that, in many interesting cases, the description of a space as a set of points is inadequate. Think for example of quantum mechanics, where position and momentum are replaced by non-commuting operators: as a consequence, Heisenberg uncertainty relations impose limitations in the precision of their simultaneous measurement so that the notion of ``point in phase space'' loses any operational meaning. Taking into account General Rela\-tivity, one can show by simple arguments that not only phase-space coordinates but also space-time coordinates should be non-commutative (cf.~\cite{DFR94,DFR95} and references therein) making the concept of ``points in space-time'' also problematic. Noncommutative Geometry provides ef\/f\/icient tools to study these ``spaces'' that are no longer described by a commutative algebra of coordinate functions, but by some noncommutative operator algebra~$\A$.

Losing the notion of points, one also loses the notion of distance between points. However one can still def\/ine a distance between \emph{states} of the algebra $\A$, for example with the help of a generalized Dirac opera\-tor~$D$. The latter is the starting point of Connes theory of \emph{spectral triples}~\cite{Con95}, which is the datum $(\A,\HH,D)$ of an involutive (non necessarily commutative) algeb\-ra~$\A$, a representation $\pi$ of $\A$ as bounded operators on a Hilbert space $\HH$ and a self-adjoint operator $D$, such that $[D,a]$ is bounded and $a(D - \lambda)^{-1}$ is compact for any $a\in\A$ and $\lambda\notin\mathrm{Sp}(D)$
(where the symbol $\pi$ is omitted). The \emph{spectral distance} between two states $\varphi_1$, $\varphi_2$ of $\A$ is def\/ined as \cite{Con89,Con96}
\begin{gather}\label{eq:5.1}
  \dc(\varphi_1,\varphi_2) \doteq \sup_{a\in\A}\big\{ \abso{\varphi_1(a) - \varphi_2(a)};\; ||[D,a]||_{\mathrm{op}}\leq 1\big\},
\end{gather}
where the norm is the operator norm coming from the representation of $\A$ on $\HH$. It is easy to check that~\eqref{eq:5.1} def\/ines a distance in a strict mathematical sense except that it may be inf\/inite.

Recall that a state is by def\/inition a positive linear application $\bar{\A}\to\C$ with norm $1$, where
$\bar{\A}$ is the $C^*$-algebra completion of $\A$. When $\A$ is the algebra of observables of a physical system, this notion of state coincides with physicist's intuition, namely density matrices or Gibbs states in statistical physics, state vectors or wave functions in quantum mechanics. For the commutative algebra $C^\infty_0(\M)$~-- i.e.~complex smooth functions vanishing at inf\/inity on some spin manifold $\M$ of dimension\footnote{Unless otherwise specif\/ied, in all the paper we assume that ``spin manifolds''  are Riemannian, f\/inite-dimensional, connected, complete, without boundary.} $m$~-- a canonical spectral triple is
\begin{gather}
  \label{eq:1}
 \big(C^\infty_0(\M),  L^2(\M,\mathcal{S}),  D\big),
\end{gather}
where $L^2(\M,\Sc)$ is the Hilbert space of square integrable spinors on which $C_0^{\infty}(\M)$ acts by point-wise multiplication
and $D$ is the Dirac operator. The latest is given in local coordinates~by
\[
 D = -i\sum_{\mu=1}^m\gamma^\mu\nabla_\mu ,
\]
where $\nabla_\mu = \partial_\mu +\omega_\mu$ is the covariant derivative associated to the spin connection $1$-form $\omega^\mu$ and $\{\gamma^\mu\}_{\mu=1,2,\dots ,m}$ are the self-adjoint Dirac gamma matrices satisfying
\begin{gather}
  \gamma^\mu\gamma^\nu + \gamma^\nu\gamma^\mu = 2g^{\mu\nu}{\mathbb I}\label{eq:67}
\end{gather}
with $g$ the Riemannian metric of $\M$ and ${\mathbb I}$ the identity matrix of dimension $2^n$ if $m=2n$ or $2n+1$.
In this case $\bar{\A} = C_0(\M)$ and the spectral distance between pure states (i.e.~states that cannot be written as a convex sum of other states), which by Gelfand transform are in one-to-one correspondence with the points of the manifold
\begin{gather}
x\mapsto\delta_x:\ \delta_x(f)\doteq f(x) \qquad \forall \; f\in C_0(\M),\label{eq:16}
\end{gather}
coincides \cite{CL92} with the geodesic distance $d$ of $\M$,
\begin{gather}\label{eq:9}
\dc(\delta_x,\delta_y)=d(x,y) \qquad \forall\;x,y\in\M.
\end{gather}

In a noncommutative framework it is tempting~-- inspired by (\ref{eq:16})~-- to take the set $\pa$ of pure states of $\A$ as the noncommutative analogue of points and $d_D$ as natural generalization of the geodesic distance. This idea has been tested in several examples inspired by physics:  f\/inite-dimensional algebras \cite{bls,IKM01}, functions on $\M$ with value in a matrix algebra \cite{MW02,Mar06,Mar08}  encoding the inner structure of the space-time of  the standard model of particles physics  \cite{Chamseddine:2007oz}, non-commutative deformations of $C^\infty_0(\M)$~\cite{Moyal}. Most often the computation of the supremum in~(\ref{eq:5.1}) is quite involved, however several explicit results have been obtained.
They all indicate that as soon as~$\A$ is noncommutative one looses an important feature
of the geodesic distance, namely~$\pa$ is no longer a path metric space~\cite{Gromov:1999fk}. Explicitly there is no curve $[0,1]\ni t\mapsto \varphi_t\in \pa$ such that
\begin{gather}
  \label{eq:10}
  d_D (\varphi_s , \varphi_t) = |t-s| d_D (\varphi_0 , \varphi_1),
\end{gather}
not even a sequence of such curves $\varphi^n$ such that
\[
  d_D (\varphi_0 , \varphi_1) = \underset{n}{\text{Inf}}  \left\{ \text{length of } \varphi^n \text{ between }\varphi_0  \text{ and }\varphi_1\right\}.
\]
This can be seen on a simple noncommutative examples studied in \cite{IKM01}, based on $\A = M_2(\C)$ acting on $\C^2$ with $D=D^*\in M_2(\C)$. $\mathcal{P}(M_2(\C))$ is weak* homeomorphic to
the Euclidean $2$-sphere. The image under this homeomorphism of the pure states $\omega_i = (\psi_i, . \psi_i)$, $i=1,2,$ def\/ined by the eigenvectors $\psi_i$ of $D$ are antipodal,  and determine a distinguished vertical axis on $S^2$. The spectral distance is invariant by rotation around this
axis and the connected components (see \eqref{eq:12} below)
are circles parallel to the horizontal plane. Explicitly, on the circle with radius $r\in[0,1]$ one f\/inds
\begin{gather}
d_D(\theta_1,\theta_2) =
    \frac{2r}{|D_1-D_2|}\left|\,\sin\frac{\theta_1 -\theta_2}{2}\right|,
\label{eq:14}
\end{gather}
where  $\theta\in[0,2\pi[$ is the azimuth and $D_i$ are the eigenvalues of $D$.
One then checks\footnote{A curve $t\mapsto\theta(t)$ satisf\/ies \eqref{eq:10} if and only if
\begin{gather}\label{eq:controesempio}
\bigl|\sin\tfrac{\theta(t) - \theta(s)}2\bigr| =K|t-s|
\end{gather}
for any $t,s\in [0,1]$ and $K$ a~constant. The right hand side of \eqref{eq:controesempio} being a~function of $t-s$, there exists a~func\-tion~$f$ such that the left hand side is $f(t-s)$. Putting $s=0$, one obtains
$f(t) = \abso{\sin \frac{\theta(t) - \theta(0)}2}$. Reinserted in~\eqref{eq:controesempio}, this yields Cauchy's functional equation $\theta(t)-\theta(s)=\theta(t-s) - \theta(0)$,
whose continuous solutions are linear. Since~\eqref{eq:controesempio} has no linear solutions, this proves the claim.}
that \eqref{eq:10} has no solution in $\mathcal{P}(M_2(\C))$~\cite{Mardev}. The lack of geodesic curves \eqref{eq:10} within $\pa$ is cured by consi\-de\-ring non-pure states. Indeed~\eqref{eq:5.1} not only generalizes the geodesic distance to noncommutative algebras, it also extends the distance to objects that are not equivalent to points, namely non-pure states. Noticing that (\ref{eq:14}) is the geodesic distance within the Euclidean disk of radius $\frac{2r}{\abso{D_1-D_2}}$ and that the set $\mathcal{S}(M_2(\C))$ of states of $M_2(\C)$ is homeomorphic to the 2-dimensional Euclidean ball (see Section~\ref{moyalsection}), one easily obtains a curve in $\mathcal{S}(M_2(\C))$ satisfying~\eqref{eq:10}, namely
\[
\varphi_t = (1-t) \varphi_0 + t\varphi_1.
\]
This remains true in full generality since, whatever algebra $\A$ and operator $D$,
\begin{gather}\label{eq:15}
d_D (\varphi_s, \varphi_t)
=\sup_{a\in\A}\big\{\abso{ (s-t)(\varphi_0  -\varphi_1 )(a)};\, ||[D,a]||_{\mathrm{op}}\leq 1\big\}
=|s-t|d_D(\varphi_0, \varphi_1).
\end{gather}

In other terms, in view of \eqref{eq:15} and \eqref{eq:9}, the spectral distance is a natural gene\-ralization of the geodesic distance to the noncommutative setting as soon as one takes into account the whole space of states, and not only its extremal points.
This motivates the study of the spectral distance between non-pure states that we undertake in this paper. We begin by giving a detailed proof of Rief\/fel's remark~\cite{Rie99} (also mentioned in~\cite{BV01}), according to which in the commutative case $\A=C_0^\infty(\M)$ the spectral distance coincides with a distance well known in optimal transport theory, namely the Wasserstein distance $W$ of order~$1$ (see the bibliographic notes of \cite[Chapter~6]{Vil08}). We stress in particular that  $\M$ has to be complete for that result to hold besides the compact case. Then we present few explicit calculations of $d_D$ between non-pure states, and question on simple examples~-- including the standard model~--  the pertinence of the optimal transport interpretation of the spectral distance in a noncommutative framework.  These are preliminaries results, intending to bring the attention of the
transport theory community on the metric aspect of noncommutative geometry, and vice-versa.

Notice that some properties of the spectral distance between non-pure states have been investigated in \cite{Rie99}: considering instead of $||[D,a]||_{\mathrm{op}}$ an arbitrary semi-norm on $\A$, it is shown that the knowledge of the distance between pure states of a noncommutative $\A$ may not be enough to recover the semi-norm on $\A$. One also needs the distance between non-pure states. This suggests that the metric information encoded in formula~\eqref{eq:5.1} is not exhausted once one knows the distance between pure states.

The plan of the paper is the following. In Section~\ref{sec:2} we recall some basics of transport theory and noncommutative geometry in order to establish~-- in Proposition~\ref{prop1}~-- the equality between $\dc$ and $W$ for any complete spin manifold $\M$. We also discuss various def\/initions of the spectral distance, characterize its connected components and emphasize the importance of the completeness condition. In Section~\ref{sec:3} we provide some lower and upper bounds for the distance.  Specializing to $\M=\R^n$, we explicitly compute the distance between a class of states generalizing Gaussian wave-packets. Section~\ref{sec:4} deals with noncommutative examples. We show that on the truncations of the Moyal plane introduced in~\cite{Moyal}, the Wasserstein distance $W_D$  with cost $d_D$ def\/ined on $\pa$ does not coincide with the spectral distance on $\sa$. But on almost-commutative geometries, including the standard model of elementary particles, the spectral distance between certain classes of states may be recovered as a Wasserstein distance $W'$ with cost~$d'$ def\/ined on a subset of $\sa$ containing $\pa$. We also point out a reformulation of the spectral distance on $\pa$ in term of the minimal work $W_I$ associated to a cost $c_I$ non-vanishing on the diagonal (i.e.~$c_I(x,x)\neq 0$).

\section{Spectral distance as Wasserstein distance of order 1}\label{sec:2}

\subsection{Spectral distance from Kantorovich duality}

For any locally compact Hausdorf\/f topological space $\X$, states $\varphi\in \Sc(C_0({\X}))$ are given by Borel probability measures $\mu$ on ${\mathcal X}$, via the rule
\begin{gather}\label{varphi}
  \varphi(f) \doteq \int_{{\mathcal X}} f \de\mu\qquad \forall \; f\in \A .
\end{gather}
This is a simple application of Riesz representation and Hahn--Banach theorems together with the assumption that ${\mathcal X}$ is $\sigma$-compact in order to avoid regularity problems. Any such $\mu$ def\/ines a state since $f$ vanishing at inf\/inity (hence being bounded) guarantees that~\eqref{varphi} is f\/inite. Pure states correspond to Dirac-delta measures.
To provide some physical intuition, one can view~$\varphi$ as a wave-packet and imagine that it describes the probability distribution of a bunch of particles. Strictly speaking a wave-packet is a square root of the Radon--Nikodym derivative of~$\de\mu$ with respect to some f\/ixed $\sigma$-f\/inite positive measure $\de x$ on $\mathcal{X}$ (it is a square-integrable function, almost everywhere def\/ined and unique modulo a phase, whenever $\de\mu$ is absolutely continuous with respect to $\de x$). For instance in quantum mechanics, with ${\mathcal X}=\R^n$ and  $\de x$ the Lebesgue measure, a wave-packet is a function $\phi\in L^2(\R^n)$ and the corresponding measure is~$\de\mu = |\phi(x)|^2\de x$.

Assuming ${\mathcal{X}}$ is a metric space with distance function $d$, there is a natural way to measure how much two states $\varphi_1$
and $\varphi_2$ dif\/fer, which is the expectation value of the distance between the two corresponding distributions
\[
\mathbb{E}(d;\mu_1\times\mu_2)=
\int_{\mathcal{X}\times\mathcal{X}}d(x,y) \de\mu_1(x)\de\mu_2(y) ,
\]
where $\mu_i$ is associated to $\varphi_i$ via \eqref{varphi}.
Other ways are suggested by transport theory (all material here is taken from \cite{Amb00,Vil03,Vil08,Brenier}, where an extensive bibliography can be found. Following \cite{Vil03} we assume from now on that $\X$ has a countable basis, so that it is a Polish space). Assume there exists a positive real function $c(x,y)$~-- the ``cost function''~-- that represents the work needed to move from $x$ to $y$. A good measure on how much the $\varphi_i$'s dif\/fer is given by the minimal work~$W$ required to move the bunch of particles from the conf\/iguration~$\varphi_1$ to the conf\/iguration~$\varphi_2$, namely
\begin{gather}
W(\varphi_1,\varphi_2)\doteq
\inf _{\pi} \int_{\mathcal{X}\times\mathcal{X}} c(x, y) \de\pi,
\label{eq:7}
\end{gather}
where the inf\/imum is over all measures $\pi$ on ${\mathcal{X}}\times {\mathcal{X}}$ with marginals $\mu_1$, $\mu_2$ (i.e.~the push-forwards of $\pi$ through the projections $\mathbb{X},\mathbb{Y}:\mathcal{X}\times\mathcal{X}\to\mathcal{X}$, $\mathbb{X}(x,y)\doteq x$, $\mathbb{Y}(x,y)\doteq y$, are $\mathbb{X}_*(\pi)= \mu_1$ and $\mathbb{Y}_*(\pi)=\mu_2$). Such measures are called \textit{transportation plans}. Finding the optimal transportation plan, that is the one which minimizes $W$, is a non-trivial question known as the Monge--Kantorovich problem. This is a generalization of Monge~\cite{Mon81} ``d\'eblais et remblais'' problem, where one considers only those transportation plans that are supported on the graph of a \emph{transportation map}, i.e.~a map $T:\mathcal{X}\to\mathcal{X}$ such that $T_{*}\mu_1=\mu_2$. Namely,
\begin{gather}
W_{\text{Monge}}(\varphi_1,\varphi_2)\doteq \inf _{T} \int_{\mathcal{X}} c(x, T(x)) \de \mu_1(x).
\label{eq:f}
\end{gather}
One of the interests of Kantorovich's generalization~\cite{Kan42} is that the inf\/imum in (\ref{eq:f}) is not always a minimum: an optimal transportation map may not exist. On the contrary the inf\/imum in (\ref{eq:7}) is a minimum and always coincides with Monge inf\/imum, even when the optimal transportation
map does not exist. Moreover when the cost function $c$ is a distance $d$,  (\ref{eq:7})  is in fact a distance on the space of states~-- with the inf\/inite value allowed~-- called the Kantorovich--Rubinstein distance (this case was f\/irst studied in~\cite{KR58}). To be sure it remains f\/inite (see \eqref{eq:3} below), it is convenient to restrict to the set  $\Sc_1(C_0(\X))$ of states whose moment of order $1$ is f\/inite, that is those distributions $\mu$ such that
\begin{gather}
\mathbb{E}(d(x_0, \circ);\mu) = \int_{\X} d(x_0,x) \de\mu(x) < +\infty,
\label{eq:2}
\end{gather}
where $x_0$ is an arbitrary but f\/ixed point in $\X$. Note that as soon as $\mathbb{E}(d(x_0,\circ);\mu)$ is f\/inite for~$x_0$, then by the triangle inequality $\mathbb{E}(d(x,\circ);\mu)\leq\mathbb{E}(d(x_0,\circ);\mu)+d(x_0,x)$ is f\/inite for any $x\in\X$ so that $S_1(C_0(\X))$ is independent on the choice of $x_0$. The Kantorovich--Rubinstein distance is also known as the \emph{Wasserstein distance of order $1$}.
The distance of order $p$ is given by a similar formula with $\mathbb{E}(d^p ;\mu_1\otimes\mu_2){^{1/p}}$, $1\leq p<\infty$, but in this paper
we are interested only in the distance of order one and we shall simply call $W$ the Wasserstein distance.

The link with non-commutative geometry, which seems to have been f\/irst noted for $\M$ compact in \cite{Rie99}, is the following: when $\X$ is a spin manifold $\M$, the Wasserstein distance with cost function the geodesic distance $d$ is nothing but the spectral distance~(\ref{eq:5.1}) associated to the canonical spectral triple~(\ref{eq:1}). This is \emph{a priori} not dif\/f\/icult to see: On one side a central result of transport theory, Kantorovich duality, provides a dual formulation of the Wasserstein distance as a supremum instead of an inf\/imum, namely (cf.~Theorem~5.10 and equation~(5.11) of~\cite{Vil08})
\begin{gather}\label{eq:8}
  W(\varphi_1, \varphi_2) = \sup_{\norm{f}_{\mathrm{Lip}}\leq 1, \,f\in L ^1(\mu_1)\cap L ^1(\mu_2) }\left(\int_\X f\de\mu_1 - \int_\X f\de\mu_2\right)
\end{gather}
for any pair of states in $S(C_0(\X))$ such that the right-hand side in the above expression is f\/inite.
The supremum is on all real $\mu_{i=1,2}$-integrable functions $f$ that are 1-Lipschitz, that is to say
\[
  \abso{f(x) -f(y)} \leq d(x,y) \qquad\forall\; x,y \in \X.
\]
On the other side, the commutator $-i[\gamma^\mu \partial_\mu, f]$ (where we use Einstein summation convention and sum
over repeated indices) acts on $L^2(\M, \mathcal{S})$ as multiplication by $-i\gamma^\mu\partial_\mu f$. Moreover  the supremum in~\eqref{eq:5.1} can be searched on self-adjoint elements \cite{IKM01}, that for $\A=\coi$ simply means real functions $f$. Thus
\begin{gather}
\norm{[D,f]}_{\text{op}}^2 =\norm{\gamma^\mu\partial_\mu f}^2_{\text{op}}=\norm{(\gamma^\mu\partial_\mu f)(\gamma^\nu\partial_\nu f)}_{\text{op}}
=\norm{\tfrac 12 (\gamma^\mu\gamma^\nu + \gamma^\nu\gamma^\mu) \partial_\mu f \partial_\nu f}_{\text{op}} \notag\\
\phantom{\norm{[D,f]}_{\text{op}}^2 }{}
=\norm{ \mathbb{I} g^{\mu\nu} \partial_\mu f \partial_\nu f}_{\text{op}}
 =\norm{ g^{\mu\nu} \partial_\mu f \partial_\nu f}_{\infty} = \norm{f}^2_{\text{Lip}}, \label{eq35ter}
\end{gather}
where we used: the $C^*$-property of the norm,  that $\partial_\mu f$ commutes with $\partial_\nu f$ and the $\gamma$'s matrices, equation~\eqref{eq:67}, that $\sqrt{ g^{\mu\nu} \partial_\mu f \partial_\nu f}$ evaluated at $x$ is the norm of the gradient $\nabla f(x)$ and
$\sup_{x\in\M}\norm{\nabla f(x)}_{T_xM} = \norm{f}_{\mathrm{Lip}}$
by Cauchy's mean value theorem. Consequently the com\-mu\-ta\-tor-norm condition in the spectral distance formula yields on $f$ the condition required in Kantorovich's
dual formula. However one has to be careful that, although the $\varphi_i$'s on the l.h.s.\ of~(\ref{eq:8}) denote states of $\coi$, the supremum on the r.h.s.\ includes functions non-vanishing at inf\/inity. Therefore~(\ref{eq:5.1}) equals~(\ref{eq:8}) if and only if the supremum on $1$-Lipschitz smooth functions vanishing at inf\/inity
is the same as the supremum on $1$-Lipschitz continuous functions non-necessarily vanishing at inf\/inity. That $C^\infty_0(\M)$ is dense within $C_0(\M)$ is well known, but it might be less known that continuous $K$-Lipschitz functions
can be approximated by smooth $K$-Lipschitz functions.
In the following we use this result -- proved for f\/inite-dimensional manifolds in \cite{GW79} and extended to inf\/inite dimension
in~\cite{AFLR07}~--  in order to to prove Rief\/fel's remark, generalized to complete locally compact manifolds (e.g.\ complete f\/inite-dimensional manifolds).

\begin{proposition}\label{prop1}
For any $\varphi_1, \varphi_2 \in \Sc(C^\infty_0(\M))$  with $\M$ a $($complete, Riemannian, finite dimensional, connected, without boundary$)$ spin manifold, one has
\[
  W(\varphi_1,\varphi_2) = \dc(\varphi_1,\varphi_2)  .
\]
\end{proposition}
\begin{proof}
i) It is well known \cite{AFLR07} that on $\R^n$ $K$-Lipschitz functions are the uniform limit of smooth
$K$-Lipschitz functions (for any $K\geq 0$).
It may be less known that the same is true for any (f\/inite-dimensional) Riemannian manifold $\M$,
and for $K$-Lipschitz functions vanishing at inf\/inity. Let us give a short proof of this result.

According to Theorem 1 in \cite{AFLR07} (that is valid for separable Riemannian mani\-folds,
so in particular for f\/inite-dimensional manifolds), given a Lipschitz function $f$,
for any continuous function $\varepsilon:\M\to\R^+$ and for any $r>0$ there exists a smooth function
$g_{\varepsilon,r}:\M\to\R$ such that
\[
 \norm{g_{\varepsilon,r}}_{\mathrm{Lip}}
\leq \norm{f}_{\mathrm{Lip}}+r \qquad \text{and}\qquad
 |f(x)-g_{\varepsilon,r}(x)|\leq\varepsilon(x) \qquad \forall \; x\in\M.
\]
As a corollary, if $f$ is a $K$-Lipschitz function vanishing at inf\/inity,
we can f\/ix a sequence $\varepsilon_n$ of continuous functions vanishing
at inf\/inity and uniformly converging to zero, and a sequence of positive
numbers $r_n$ converging to zero in order to get a sequence of smooth functions $g_n:\M\to\R$ such that
\[
 \norm{g_n}_{\mathrm{Lip}}
\leq K+r_n \qquad \text{and}\qquad |f(x)-g_n(x)|\leq\varepsilon_n(x) \qquad \forall \; x\in\M .
\]
Obviously $g_n$ vanishes at inf\/inity
(since both $f$ and $\varepsilon_n$ vanish at inf\/inity). Let us call
$f_n:=K(K+r_n)^{-1}g_n$. The $f_n$'s are the required smooth $K$-Lipschitz functions
vanishing at inf\/inity that converge uniformly to $f$.
Indeed
\[
f-f_n=\frac{K}{K+r_n}(f-g_n)+\frac{r_n}{K+r_n}f
\]
implies
\[
\sup_x |f(x)-f_n(x)|\leq \frac{K}{K+r_n}\sup_x |\varepsilon_n(x)|+\frac{r_n}{K+r_n}\sup_x|f(x)| ,
\]
and the right hand side goes to zero for $n\to\infty$ since  $\sup_x |\varepsilon_n(x)|\to 0$
and $r_n\to 0$ by assumption, while $\sup_x|f(x)|$ is f\/inite since $f$ is continuous
vanishing at inf\/inity.

ii) By (\ref{eq35ter}) the supremum in (\ref{eq:5.1}) is on $1$-Lipschitz smooth functions vanishing at inf\/inity. By i) above any $1$-Lipschitz function $f$ in $C_0(\M)$ can be uniformly approximated by smooth $1$-Lipschitz functions $f_n$ in $C^\infty_0(\M)$. Since any state~$\varphi$
of $C_0(\M)$ is continuous \cite{BR} with respect to the sup-norm, namely
\[
  \lim_{n\rightarrow +\infty}
\norm{f - f_n }_{\infty} = 0\quad \Longrightarrow \quad
  \lim_{n\rightarrow +\infty}\varphi(f_n) = \varphi(f),
\]
the supremum in (\ref{eq:5.1}) can be equivalently searched on continuous functions and the spectral distance writes
\begin{gather}
 \dc(\varphi_1, \varphi_2) = \sup_{f\in C_0(\M,\R) ,\,\norm{f}_{\mathrm{Lip}}\leq 1}
   \left(\int_\M f\de\mu_1 - \int_\M f\de\mu_2\right)  .\label{eq:4}
\end{gather}

iii) In case $\M$ is compact $C_0(\M,\R)$ coincides with $ C(\M,\R)$ and
\eqref{eq:4} equals \eqref{eq:8}, hence the result. In case $\M$  is only locally compact, $C_0(\M,\R)\subset C(\M,\R)$ so that
$d_D(\varphi_1, \varphi_2 ) \leq W(\varphi_1, \varphi_2 )$.
To get the equality, it is suf\/f\/icient to show that to any $1$-Lipschitz $\mu_i$-integrable function $f\in  C(\M,\R)$ is associated a sequence
of functions $f_n \in C_0(\M,\R)$ such that
\begin{gather}
||f_n||_{\mathrm{Lip}}\leq 1
\label{eq:23}
\end{gather}
and
\begin{gather}
\label{eq:24}
\lim_{n\rightarrow +\infty}  \Delta(f_n) = \Delta(f),
\end{gather}
where  $ \Delta(f)\doteq \varphi_1(f) -\varphi_2(f)$.
We claim that such a sequence is given by
\begin{gather}
f_n(x)\doteq f(x)e^{-d(x_0,x)/n}  ,\qquad n\in\N,
\label{eq:010}
\end{gather}
where $x_0$ is any f\/ixed point.
Indeed since $\M$ is complete, $d(x_0,x)$ diverges at inf\/inity as explained in Section~\ref{completesection} below;
by the $1$-Lipschitz condition $|f(x)|\leq |f(x_0)|+d(x_0,x)$, and this proves that
$|f_n(x)|\leq \bigl(|f(x_0)|+d(x_0,x)\bigr)e^{-d(x_0,x)/n}$ vanishes at inf\/inity.

To obtain \eqref{eq:23}, one f\/irst notices
that $\Delta(f+C)=\Delta(f)$ for any $C\in\R$, so that we can
assume without loss of generality that $f(x_0)=0$, that is to say
\[
|f(x)| = |f(x)-f(x_0)|\leq d(x_0,x)  .
\]
Second, from
$\nabla f_n=\left(\nabla f-n^{-1}f  \nabla  d(x_0,\circ\,)\right)
e^{-d(x_0, \circ )/n}$ and remembering that both $f$ and
$d(x_0, \circ )$ are $1$-Lipschitz, one gets
\[
|\nabla f_n| \leq\big(1+n^{-1} d(x_0,\circ )\big)e^{-n^{-1} d(x_0,\circ)}  .
\]
The inequality
\eqref{eq:23}  then follows from $(1+\xi)e^{-\xi}\leq 1$ $\forall\; \xi>0$.
\eqref{eq:24} comes from Lebesgue's dominated convergence theorem:
$|f_n(x)|\leq |f(x)|$ $\forall\; x, \, n$ and $|f|$ is $\mu_i$ integrable by hypothesis, so that
$\lim\limits_{n\rightarrow \infty} \int_\M f_n \de\mu_i = \int_\M f \de\mu_i$.
\end{proof}

\subsection{Alternative def\/initions}

There exist several equivalent formulations of the spectral distance. First one may consider
continuous instead of smooth functions. Indeed, as explained in \cite{Con94}, for any measurable bounded function $f$
one can view $[D,f]$ as the bilinear form
\[
  \xi, \eta \rightarrow \langle D\xi, f\eta \rangle - \langle f^*\xi, D\eta \rangle  ,
\]
 well def\/ined on the domain of $D$ (a dense subset of $L^2(\M,\mathcal{S})$). Therefore $[D,f]$ makes sense also when
$f$ is not smooth and one can def\/ine~\cite{Con95} the spectral distance as the supremum on all continuous functions $f\in C_0(M)$ with $\norm{[D,f]}\leq 1$, that is the set of $1$-Lipschitz functions, obtaining thus directly~\eqref{eq:4}.
In the literature one f\/inds both def\/initions: supremum on continuous functions~\cite{Con95,fgbv} or on smooth functions~\cite{Con96,Connes:2008kx}.

Second, one may be puzzled by the use of the spin structure to recover the Riemannian metric. In fact instead of the Dirac operator one can equivalently use, as explained in \cite{connesckm}, the signature ope\-rator $d+d^\dag$ acting on the Hilbert space $L^2(\M, \wedge)$ of square-integrable dif\/ferential forms, or the de-Rham Laplacian $\Delta = dd^\dag + d^\dag d$ acting on $L^2(\M)$.
Here $d$ is the exterior derivative and~$d^\dag$ its adjoint with respect to the inner product~\cite{fgbv}
\[
  \langle\omega, \omega' \rangle = \int_\M  (\omega, \omega') \nu_g \qquad\forall\; \omega, \omega'\in  L^2(\M, \wedge)
\]
with $\nu_g$ the volume form  and the inner product on $k$-form given by
\[
\big(dx^{\alpha_1}\wedge \dots  \wedge dx^{\alpha_k}, dx^{\beta_1}\wedge \dots  \wedge dx^{\beta_{k'}}\big) = \delta_{kk'}\det \big(g^{\alpha_i\beta_j}\big) \qquad
1\leq i,j\leq k.
\]
The action of both operators only depends on the Riemannian structure and suitable commutators yield on self-adjoint elements $f=f^*$ in $C^\infty_0(\M)$ the same semi-norm as the commutator with $D$. Explicitly\footnote{Although we use the same symbol, one should keep in mind that the three operator norms in \eqref{various} correspond to actions on dif\/ferent Hilbert spaces: $L^2(\M, \wedge)$, $L^2(\M)$, $L^2(\M,\mathcal{S})$.}
\begin{gather}
\label{various}
\norm{[d+d^\dag,f]^2}_{\mathrm{op}} = \tfrac{1}{2} \norm{[[\Delta ,f] ,f]}_{\mathrm{op}} =\norm{[D,f]}_{\mathrm{op}}^2.
\end{gather}
To show these equalities,  let us note that on $L^2(\M, \wedge)$,
$[d,f] = \epsilon(df)$ where $\epsilon$ denotes the wedge multiplication,
\[
  \epsilon(df)\omega \doteq df\wedge \omega.
\]
Therefore $[d,f]$  commutes with $0$-form, and the same is true for its adjoint $[d,f]^\dag$. Assuming $f=f^*$, that is $[d^\dag,f] = -[d,f]^\dag$,
few manipulations with commutators yield
\begin{subequations}
\begin{gather}
  [[\Delta,f],f]  = [[d,f]d^\dag,f] + [d[d^\dag,f],f] +[[d^\dag,f]d,f] + [d^\dag[d,f],f] \\
\phantom{[[\Delta,f],f]}{}
= 2 [d,f][d^\dag,f] + 2 [d^\dag,f][d,f]  \label{eq:35}\\
\phantom{[[\Delta,f],f]}{}
=2[d+d^\dag,f]^2, \label{eq:35bis}
\end{gather}
\end{subequations}
where  $[d,f]^2 = [d^\dag,f]^2=0$ due to the graded commutativity of the wedge product. Remembering that $\epsilon(dx^\mu)^\dag = g^{\mu\nu}\iota(\partial_\nu)$ where $\iota$ is the contraction (see e.g.~\cite{fgbv}), one gets
$
  [d,f]^\dag = \epsilon(df)^\dag = ((\partial_\mu f) \epsilon(dx^\mu))^\dag = (g^{\mu\nu}\partial_\nu f) \iota(\partial_\mu)
$
so that~(\ref{eq:35}) becomes
\begin{gather*}
  [[\Delta,f],f]  = -2 \partial_\rho f  (g^{\mu\nu}\partial_\nu f) \left(\epsilon(dx^\rho)\iota(\partial_\mu) + \iota(\partial_\mu) \epsilon(dx^\rho)\right) =
-2g^{\nu\rho} \partial_\rho f \partial_\nu f,
\end{gather*}
where we used the fact that $\epsilon(dx^\rho)\iota(\partial_\mu) + \iota(\partial_\mu) \epsilon(dx^\rho) = \delta^{\rho}_\mu$.
With (\ref{eq:35bis}) this shows that  $[d+d^\dag,f]^2 = \frac 12 [[\Delta ,f] ,f]$ is the operator of point-wise
multiplication by the function $g^{\mu\nu}\partial_\mu f\partial_\nu f$, whose operator norm is $\norm{g^{\mu\nu} \partial_\mu f\partial_\nu f}_{\infty}$. (\ref{various}) then follows from (\ref{eq35ter}).

Notice that instead of the Laplacian one could use any other 2nd order dif\/ferential operator with
the same principal symbol.

\subsection{On the importance of being complete}\label{completesection}

At point iii~in the proof of Proposition \ref{prop1} it is crucial that $\M$ is complete.  By the Hopf--Rinow theorem a complete f\/inite-dimensional  Riemannian manifold  is a proper metric space, hence the geodesic distance from any f\/ixed
point $x_0$, $x\mapsto d(x_0,x)$, is a proper map~\cite{Roe96}. In particular for non-compact $\M$ this means
that $d(x_0,x)$ diverges at inf\/inity, so that the functions $f_n$ in \eqref{eq:010}
vanish at inf\/inity.

When $\M$ is not complete, not only the $f_n$ do not vanish at inf\/inity which spoils the proof, but also the
def\/inition of the Wasserstein distance requires more attention. Indeed in \cite{Vil03} Kantorovich duality is proved assuming $\X$ is complete. It is not clear to the authors whether the duality holds in the non-complete case. However one can still take (\ref{eq:8})  as a def\/inition of~$W$, letting aside whether this is the same quantity as (\ref{eq:7}) or not. Then it is easy to see that on non-complete~$\M$ the spectral distance and~$W$ are not necessarily equal. Suppose indeed that $\NN$ is compact, and $\M$ is obtained from $\NN$ by
removing a point~$x_0$, so that~$\M$ is locally compact and \emph{not} complete. The spectral distances $d_D^\M$ and  $d_D^\NN$ computed on $\M$ and $\NN$ are equal.
Indeed by~\eqref{eq:4}, which holds also in the non-complete case,
we can compute $\dc$ as supremum over continuous functions instead of smooth ones;
in the computation of the spectral distance on~$\NN$, it is equivalent to take the supremum
over $1$-Lipschitz $f\in C(\NN)$ or over $1$-Lipschitz $f'=f-f(x_0)\in C(\NN)$ vanishing at $x_0$
(since $\varphi_1(f)-\varphi_2(f)=\varphi_1(f')-\varphi_2(f')$ for any two states $\varphi_1,\varphi_2$);
therefore
\begin{gather*}
  d^{\NN}_D(\varphi_1,\varphi_2)  = \sup_{f\in C(\NN)} \big\{
  \varphi_1(f) - \varphi_2(f);\, ||f||_{\mathrm{Lip}}\leq
  1\big\}   \\
  \phantom{d^{\NN}_D(\varphi_1,\varphi_2)}{}
   = \sup_{f\in
    C(\NN),f(x_0)=0} \big\{ \varphi_1(f) -
    \varphi_2(f);\, ||f||_{\mathrm{Lip}}\leq 1\big\}  = d^{\M}_D(\varphi_1,\varphi_2),
\end{gather*}
where in last equality we used $C_0(\M)=\{f\in C(\NN),f(x_0)=0\}$.

On the contrary,  $W$ computed between pure states coincide with the geodesic distance, and the latest may or may not be the same
on $\NN$ and $\M$. For instance one does not modify the geodesic distance by removing a point from the two sphere. But taking for $\NN$  the circle, thought of as the closed interval $[0,1]$
with $0$ and $1$ identif\/ied, and $\M$ the open interval $(0,1)$,
one gets $W_{\NN}(x,y) = \min\{|x-y|,1-|x-y|\}$ whereas $W_{\M}(x,y) = |x-y|$.
Thus on $\M=(0,1)$
\[
d_D^\M(x,y) = d_D^\NN(x,y) = W_\NN(x,y)\neq W_\M(x,y).
\]

\subsection{On  the hypothesis of f\/inite moment of order 1}

Restricting to the states with f\/inite moment of order one (cf.~\eqref{eq:2}) guarantees that Wasserstein distance is f\/inite, since
by~\eqref{eq:8}
\begin{gather}
W(\varphi_1, \varphi_2)  = \sup_{\norm{f}_{\mathrm{Lip}}\leq 1}\left(\int_\X \big(f(x)-f(x_0)\big)\de\mu_1(x)
 - \int_\X \big(f(x)-f(x_0)\big)\de\mu_2(x)\right) \notag\\
\phantom{W(\varphi_1, \varphi_2)}{}
\leq\sup_{\norm{f}_{\mathrm{Lip}}\leq 1}\int_\X \big|f(x)-f(x_0)\big|\de\mu_1(x)
 +\sup_{\norm{f}_{\mathrm{Lip}}\leq 1}\int_\X \big|f(x)-f(x_0)\big|\de\mu_2(x) \notag\\
\phantom{W(\varphi_1, \varphi_2)}{}
\leq \int_\X d(x,x_0)\de\mu_1(x)+\int_\X d(x,x_0)\de\mu_2(x) <\infty  .\label{eq:3}
\end{gather}
An obvious upper-bound is then obtained by choosing $\pi=\mu_1\times\mu_2$ in \eqref{eq:7}:
\begin{gather}
  \label{eq:11}
  \dc(\varphi_1, \varphi_2) = W(\varphi_1, \varphi_2) \leq \mathbb{E}(d;\mu_1\times\mu_2) .
\end{gather}
When at least one of the states is pure, this upper bound is an exact result, even outside $\Sc_1(\coi)$.

\begin{proposition}\label{cor2}
For any $x\in \M$ and any state $\varphi\in \Sc(\coi)$,
\[
\dc(\varphi,\delta_x)= W(\varphi,\delta_x) = \mathbb{E}\bigl(d(x,\circ);\mu\bigr) .
\]
\end{proposition}

\begin{proof}
By Proposition \ref{prop1},
\begin{gather*}
\dc(\varphi,\delta_x) =  W(\varphi,\delta_x)
 =\sup_{\norm{f}_{\mathrm{Lip}}\leq 1, f\in L^1(\mu)} \left( \int_\M f\de\mu -f(x)\right)  \\
 \phantom{\dc(\varphi,\delta_x) =  W(\varphi,\delta_x)}{}
 =\sup_{\norm{f}_{\mathrm{Lip}}\leq 1,  f\in L^1(\mu)} \left( \int_\M(f(y)-f(x)\bigr)\,\de\mu(y) \right) \\
 \phantom{\dc(\varphi,\delta_x) =  W(\varphi,\delta_x)}{}
\leq \sup_{\norm{f}_{\mathrm{Lip}}\leq 1,  f\in L^1(\mu)}\int_\M\abso{ f(x)-f(y)}\,\de\mu(y)\\
 \phantom{\dc(\varphi,\delta_x) =  W(\varphi,\delta_x)}{}
\leq\int_\M d(x,y)\,\de\mu(y) =\mathbb{E}\bigl(d(x,\circ);\mu\bigr).
\end{gather*}
This supremum is attained on the $1$-Lipschitz functions $f(y)\doteq d(x,y)$ in case $\mu\in\Sc_1(\coi)$, or is obtained by the sequence $f_n$ as def\/ined in~\eqref{eq:010} in case~$\mu$ has not a  f\/inite moment of order~$1$.
\end{proof}

Notice that this proposition does not rely on the f\/initeness of $W$ nor $d_D$, and makes sense since $\mathbb{E}\bigl(d(x,\circ);\mu\bigr)$ is either
inf\/inite or convergent, the integrand in (\ref{eq:2}) being a positive function. When the distributions are both localized around two points, $x$, $y$, transportation maps are simply paths from $x$ to $y$ and the minimal work coincides with the cost $c(x,y)$ to move from $x$ to $y$, i.e.~with the geodesic distance. In other words
\begin{gather}
d_D(\delta_x, \delta_y) = W(\delta_x,\delta_y) = d(x,y)
\label{eq:43}
\end{gather}
and one retrieves \eqref{eq:9}.
Proposition \ref{cor2} also  yields an alternative def\/inition of $S_1(\coi)$.
\begin{corollary}\label{cooro1}
$\varphi\in \Sc_1(\coi)$ if and only if $\varphi$ is at f\/inite spectral-Wasserstein distance from
any pure state.
\end{corollary}
\begin{proof}
By the triangle inequality, the moment of order one of $\varphi$ is f\/inite for a f\/ixed
$x=x_0$ if and only if it is f\/inite for all $x\in\M$.
\end{proof}

Let us conclude this section with some topological remarks.

\begin{definition}\rm
Given an arbitrary spectral triple $(\A,\HH,D)$, we call \cite{Mar08}
\begin{gather*}
   \text{Con}(\varphi) \doteq \{\varphi'\in S(\A),\; d_D(\varphi, \varphi')<+\infty\} .
\end{gather*}
\end{definition}

Notice that connected components in $\mathcal{S}(\A)$ for the topology induced by $d_D$ coincide with sets
of states at f\/inite distance from each other, thus justifying the name $\text{Con}(\varphi)$.
Indeed $\mathrm{Con}(\varphi)$ is path-connected since for any $\varphi_0,\varphi_1\in\mathrm{Con}(\varphi)$
the map
\[
[0,1]\ni t\mapsto
\varphi_t = (1-t) \varphi_0 + t\varphi_1
\in \mathrm{Con}(\varphi)
\]
is continuous (for all $\epsilon>0$ called $\delta_\epsilon=\epsilon/d(\varphi_0,\varphi_1)$ from~\eqref{eq:15} we get
$|t-s|<\delta_\epsilon \Rightarrow d_D (\varphi_t, \varphi_s)  <\epsilon$).
That $\mathrm{Con}(\varphi)$ is maximal~-- i.e.~there is no connected component containing it properly~-- can
be easily seen: for any $\varphi'\in\mathrm{Con}(\varphi)$ any open ball
centered at $\varphi'$ is contained in $\mathrm{Con}(\varphi)$, so that $\mathrm{Con}(\varphi)$ is open;
by the triangle inequality the same is true for the complementary set, proving that $\mathrm{Con}(\varphi)$
is also closed. Therefore, any set containing properly $\mathrm{Con}(\varphi)$ is not connected since it
contains a subset that is both open and closed.

For the Wasserstein distance, the set of states with f\/inite moment of order one
is a connected component.

\begin{corollary}
\label{corlabel}
For any $\varphi\in \Sc_1(\coi)$, $\mathrm{Con}(\varphi)= \Sc_1(\coi)$.
\end{corollary}
\begin{proof}
Let $\varphi\in \Sc_1(\coi)$. By the triangle inequality,
\[
\dc(\varphi,\varphi')\leq\dc(\varphi,\delta_x)+\dc(\delta_x,\varphi')
\]
is f\/inite by Corollary~\ref{cooro1} as soon as $\varphi'\in \Sc_1(\coi)$.
Similarly
\[
\dc(\delta_x,\varphi')\leq\dc(\delta_x,\varphi)+\dc(\varphi,\varphi')
\]
implies that $\dc(\varphi,\varphi')$ is inf\/inite when $\varphi'\notin \Sc_1(\coi)$.
\end{proof}

Notice that considering only pure states, the set
\begin{gather}\label{eq:12}
\widetilde{\text{Con}}(\varphi)\doteq \text{Con}(\varphi)\cap \pa,\qquad\varphi\in\pa
\end{gather}
is not necessarily (path)-connected in $\pa$.
In the example~\eqref{eq:14} $\widetilde{\text{Con}}(\varphi)$ is indeed a connected component of $\pa$; but in the standard model (see Section~\ref{smsection}) $\widetilde{\text{Con}}(\varphi)=\pa$ and contains two disjoint connected components $(\M,\delta_\C), (\M,\delta_{\mathbb{H}})$.

\section{Bounds for the distance and (partial) explicit results}\label{sec:3}

Having in mind that noncommutative geometry furnishes a description of the full standard model of the electro-weak and strong interactions minimally  coupled to Euclidean general relativity~\cite{Chamseddine:2007oz}, computing the spectral distance could be a way to obtain a ``picture'' of spacetime  at the scale of unif\/ication. Regarding pure states, this picture has been worked out in~\cite{MW02} and is recalled in Section~\ref{sec:4}: one  f\/inds that the connected component of $d_D$ is the disjoint union of two copies of a spin manifold $\M$, with distance between the copies coming from the Higgs f\/ield. Extending this picture to non-pure states is far from trivial since already on Euclidean space explicit computations of the Wasserstein distance are very few. This does not seem to be the most interesting issue in optimal transport, where one is rather interesting in determining the optimal plan than computing~$W$. On the contrary from our perspective computing $d_D$ is of most interest, while f\/inding the optimal plan (i.e.~the element that reaches the supremum in the distance formula) is not an aim in itself. In this section, we collect various results on the Wasserstein-spectral distance in the commutative case $\A = C_0^\infty(\M)$: upper and lower bounds for the distance on any spin manifold $\M$, and explicit result for a certain class of states in case $\M = \R^n$. Some of these results might be known from optimal transport theory, but we believe it is still interesting to present them from our perspective.

We postpone to the next section a discussion of the noncommutative case. In all this section, $d_D=W$ and to avoid repetition we simply call it  the spectral distance.

\subsection{Upper and lower bounds on any spin manifold}

On a spin manifold, a lower bound on the spectral distance between states with f\/inite moment of order $1$ is given by the distance between their mean points. For  $\M = \R^m$ the mean point of $\varphi\in\Sc_1(\coi)$ is the barycenter of $\mu$, namely
 \begin{gather}
 \bar{x}\doteq   (\bar{x}^\alpha) \qquad \text{with} \quad
\bar{x}^\alpha =  \mathbb{E}(x^\alpha ; \mu)=\int_{\R^m} x^\alpha \de\mu, \label{eq:6}
 \end{gather}
where
\[
 x^\alpha : \ \R^m \rightarrow \R, \qquad
 x\mapsto x^\alpha(x)  ,
\]
with  $\alpha=1,\dots  ,m$ denote a set of Cartesian coordinate functions on $\R^m$. For $\M$ that is not the Euclidean space, the mean point can be def\/ined through Nash embedding, that is an isometric embedding of $\M$ onto a subset $\widetilde{\M}$ of the Euclidean space $\R^n$ for some $n\geq m$ \cite{Nas54,Nas56}. We say that the embedding $N:\M\to\widetilde{\M}$ is convex if $\widetilde{\M}$ is a convex subset of $\R^n$. In that case,  the barycenter $\tilde x$ of $\tilde{\mu} \doteq N_* \mu$
is in $\widetilde{\M}$, and
\[
 \bar{x}\doteq N^{-1}  (\tilde{x})\in\M
\]
is a well def\/ined generalization of~(\ref{eq:6}) to $\M$.

Before showing in Proposition~\ref{nashprop} that the distance between mean points bounds from below the spectral distance, let us recall two properties of Nash embedding that will be useful in the following.

\begin{lemma}\label{nashlem}
Let $\M$ be a Riemannian manifold admitting a convex isometric
embedding $N:\M\to\widetilde{\M}\subset\R^n$. Then
\begin{gather}\label{eq:22}
d(x,y) = \abso{N(x) - N(y)},
\end{gather}
where $\abso{\,\cdot\,}$ is the Euclidean distance on $\R^n$.
Moreover, if $f$ is a $1$-Lipschitz function on the metric space $(\M,d)$
then $\tilde{f}\doteq f\circ N^{-1}$ is $1$-Lipschitz on $(\widetilde{\M}, \abso{\,\cdot\,} )$.
\end{lemma}

\begin{proof}
An isometry $N:\M\to\widetilde{\M}$ preserves the distance function and
sends geodesics to geodesics (cf.~e.g.~\cite[page~61]{Cha93}). If $\widetilde{\M}$ is
convex, since the metric on $\widetilde{\M}$ is the restriction of the Euclidean metric
of $\R^n$, geodesics are straight lines and the distance is the Euclidean one.
This proves \eqref{eq:22}.

As a consequence, if $f$ is a $1$-Lipschitz function on $\M$, for any $\tilde{x},\tilde{y}\in\widetilde{\M}$
we have
\[
\abso{\tilde f (\tilde x) -\tilde f(\tilde y)}=\abso{f(x) - f(y)}
  \leq d(x,y)
  =\abso{\tilde{x}-\tilde{y}},
\]
where $x=N^{-1}(\tilde{x})$ and $y=N^{-1}(\tilde{y})$.
This means that $\tilde{f}$ is $1$-Lipschitz on $\widetilde{\M}$.
\end{proof}

\begin{proposition}\label{nashprop}
Let $\M$ be a Riemannian spin manifold that admits a convex isometric embedding $\M\hookrightarrow\R^n$.
For any states $\varphi_1$, $\varphi_2$ in $S(C_0(\M))$ with mean points $\bar{x}_1$, $\bar{x}_2$,
\[
d(\bar{x}_1,\bar{x}_2)\leq\dc(\varphi_1,\varphi_2) \leq \mathbb{E}\bigl(d ;\mu_1\times\mu_2\bigl)  .
\]
\end{proposition}

\begin{proof}
\eqref{eq:11} holds true, so we need to prove only the lower bound.
Let us f\/ix a basis of the vector space $\R^n$ such that
\begin{gather}
\tilde{x}_1 - \tilde{x}_2=\bigl(\abso{\tilde{x}_1-\tilde{x}_2},0, \dots , 0\bigr) ,\label{eq:68}
\end{gather}
where $\tilde{x}_i =N(\bar{x}_i)$, $i=1,2,$ has Cartesian coordinates
\[
\int_{\widetilde{\M}} x^\beta\,\de\tilde{\mu_i}, \qquad \beta=1, \dots ,n.
\]
\eqref{eq:68} is equivalent to
\begin{gather*}
\int_{\widetilde{\M}} {x}^1 \de \tilde{\mu}_1 - \int_{\widetilde{\M}} x^1 \de \tilde{\mu}_2  = \abso{\tilde{x}_1 - \tilde{x}_2} \quad\text{and}\quad \int_{\widetilde{\M}} {x}^\beta \de \tilde{\mu}_1 - \int_{\widetilde{\M}} {x}^\beta \de \tilde{\mu}_2 =0  \quad \text{for} \quad \beta\in\left[2,n\right].
\end{gather*}
Therefore
\begin{align*}
\abso{\tilde x _1 -\tilde x _2}
=\Big|\int_{\widetilde{\M}} {x}^1 \de \tilde{\mu}_1 - \int_{\widetilde{\M}} x^1 \de \tilde{\mu}_2\Big|
\leq\sup_{\norm{\tilde f}_{\mathrm{Lip}\leq 1}} \Big| \int_{\widetilde{\M}} \tilde f\de\tilde \mu_1 -
\int_{\widetilde{\M}} \tilde f\de\tilde \mu_2\Big|
\end{align*}
where we noticed that $x^1$ is Lipschitz in $\widetilde{\M}$ with constant $1$ since
\[
{x}^1(\tilde x) - {x}^1(\tilde y) \leq
 \sqrt{({x}^1(\tilde x) - {x}^1(\tilde y))^2+  \sum\nolimits_{\alpha=2}^n \left({x}^\alpha(\tilde x) -{y}^\alpha(\tilde x)\right)^2}= \abso{\tilde x-\tilde y}  .
\]
Using $\int_{\widetilde{\M}}\tilde{f}\de\tilde{\mu}=\int_{\M}f\de\mu$ and Lemma \ref{nashlem},
\[
\abso{\tilde x _1 -\tilde x _2}\leq\sup_{\norm{f}_{\mathrm{Lip}\leq 1}} \left|\int_\M  f\de \mu_1 -
\int_{\M} f\de \mu_2\right|=\dc(\varphi_1, \varphi_2).
\]
Proposition follows from (\ref{eq:22}), i.e.~$\abso{\tilde x _1 -\tilde x _2}=d(\bar{x}_1,\bar{x}_2)$.
\end{proof}

Note that when $\M$ does not admit a convex embedding, Proposition~\ref{nashprop} still holds with
$\abso{\tilde x_1 - \tilde x_2}$ instead of $d(\bar{x}_1,\bar{x}_2)$. However this might not be the
most interesting lower bound since it involves a distance on~$\R^n$ that is not the push-forward
of the one on $\M$ (the points $\tilde{x}_i$ may \emph{not}  be in $\widetilde{\M}$, and the Euclidean
distance, even when restricted to $\widetilde{\M}$, is \emph{not} the push-forward of the geodesic
distance on~$\M$).

\subsection{Spectral distance in the Euclidean space}

In this section we study the spectral distance between states of $\Sc_1(C^\infty_0(\M))$ in the case $\M =\R^m$.
To a given density probability $\psi\in L^1(\R^m)$ with f\/inite moment of order $1$, e.g.~a Gaussian
$\psi(x)~=~\pi^{-\frac{m}{2}}e^{-|x|^2}$, one can associate
a state $\Psi_{\sigma,x}$ given by
\[
\Psi_{\sigma,x}(f)\doteq\frac{1}{\sigma^m}\int_{\R^m}f(\xi)\psi(\tfrac{\xi-x}{\sigma})\de^m\xi  ,
\]
for any $x\in\R^m$ and $\sigma\in\R^+$. This becomes the pure state (the point) $\delta_x$ in the $\sigma\to 0^+$ limit since
$\lim\limits_{\sigma\to 0^+} \psi_{\sigma, x}(f)=f(x)$. In this sense $\Psi_{\sigma,x}$ can be viewed as a ``fuzzy'' point, that is to say a wave-packet~-- characterized by a shape $\psi$ and a
 width $\sigma$~-- describing the uncertainty in the localization around the point~$x$.

The spectral distance between wave packets with the same shape is easily calculated.

\begin{proposition}
\label{propdeloc}
The distance between two states $\Psi_{\sigma,x}$ and $\Psi_{\sigma',y}$ is
\begin{gather}\label{eq:r2dist}
\dc(\Psi_{\sigma,x}, \Psi_{\sigma',y})=\int |x-y+(\sigma-\sigma')\xi| \psi(\xi)\de^m\xi .
\end{gather}
In particular, for $\sigma=\sigma'$ the distance does not depend on the shape
 $\psi$:
\[
d_D(\Psi_{\sigma,x}, \Psi_{\sigma,y})=|x-y| .
\]
\end{proposition}

\begin{figure}[t]
\centerline{\includegraphics[width=8.5cm]{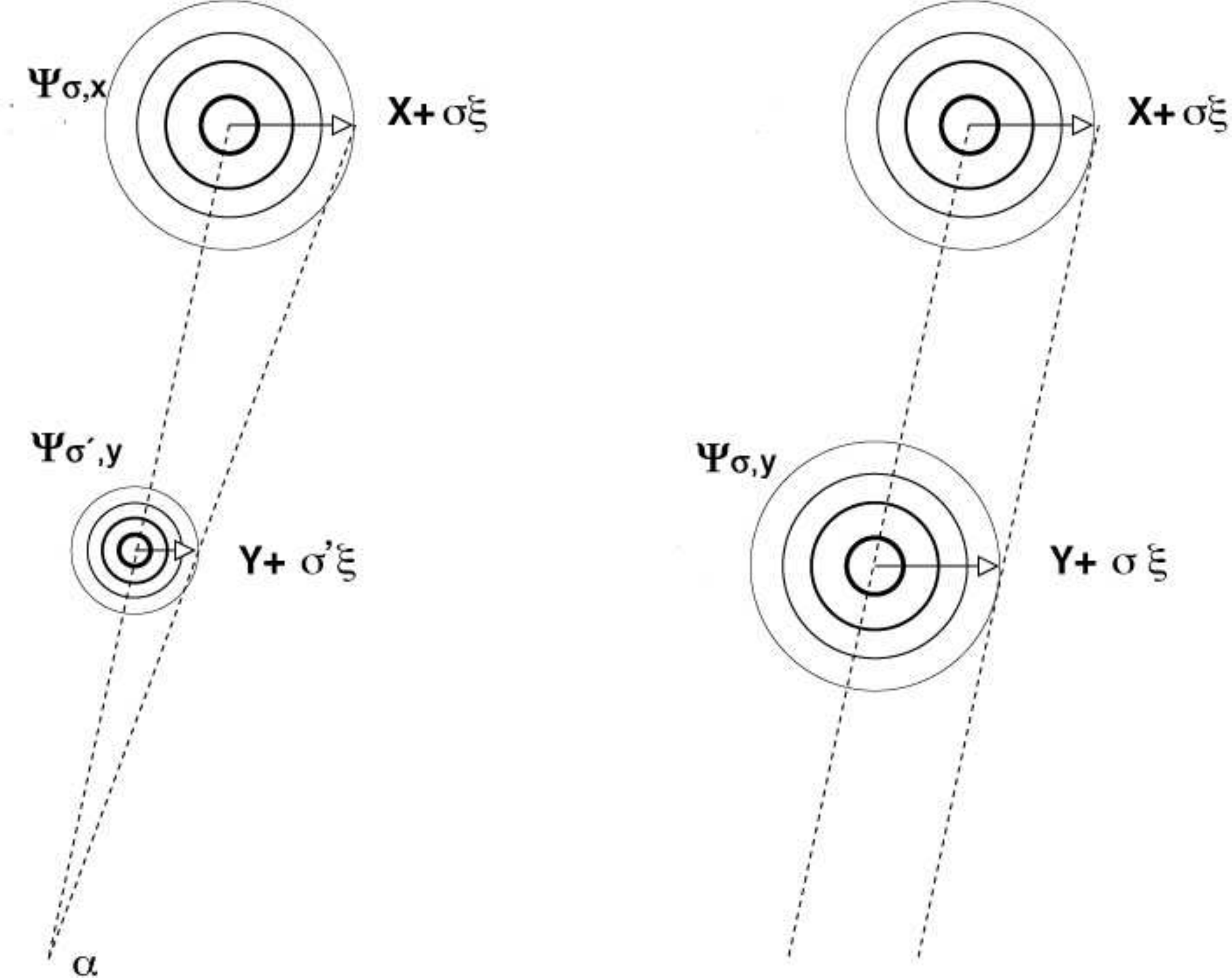}}

\caption{The function $h$ that attains the supremum, in case
$\sigma>\sigma'$ and $\sigma=\sigma'$. Dot lines are tangent to the
gradient of $h$.}\label{fig}
\end{figure}

\begin{proof}
Since $f(z)-f(w)\leq |z-w|$ for any $1$-Lipschitz function $f$, we have
\begin{gather*}
\Psi_{\sigma,x}(f)-\Psi_{\sigma',y}(f)
 =\int\left(f(\sigma\xi+x)-f(\sigma'\xi+y)\right) \psi(\xi) \de^m\xi \\
\phantom{\Psi_{\sigma,x}(f)-\Psi_{\sigma',y}(f)}{} \leq \int |x-y+(\sigma-\sigma')\xi|\,\psi(\xi)\,\de^m\xi.
\end{gather*}
When $\sigma=\sigma'$ the upper bound is attained by the function
\[
h(z)=z\cdot\frac{x-y}{|x-y|}  ,
\]
while for $\sigma\neq\sigma'$ it is attained by the function
\[
h(z)=\left|z-\alpha\right|,
\]
where
\begin{gather}\label{eq:alpha}
\alpha\doteq\frac{\sigma'x-\sigma y}{\sigma'-\sigma}
\end{gather}
as shown in Fig.~\ref{fig}.
\end{proof}

 In case $\sigma\neq\sigma'$ the function $h$ that attains the supremum
measures the geodesic distance between $z\in\R^m$ and the point
$\alpha$ def\/ined in~(\ref{eq:alpha}). Geometrically the latest is the
intersection of $(x,y)$ with $(x+\sigma\xi, y+\sigma'\xi)$, and is
independent of $\xi$. In case $\sigma=\sigma'$, $\alpha$ is send to
inf\/inity and~$h$ measures the length of the projection of~$z$ on the $(x,y)$ axis (see Fig.~\ref{fig}).
The picture is still valid for pure states,
i.e.~$\sigma=\sigma'= 0$: $h$ can be taken either as the
geodesic distance to any point on $(x,y)$ outside the segment $[x,y]$,
or as the distance to any axis perpendicular to $(x,y)$ that does not
intersect $[x,y]$.

For Gaussian shape, Proposition~\ref{propdeloc} can be confronted with
the Wasserstein distance of order~$2$ between Gaussians computed in~\cite{GS84}.
The f\/inal expression for the distance of order $2$ is quite simpler, since it is
an algebraic function of the means and covariance matrices. Here, computing
the integral in~\eqref{eq:r2dist} for $\psi$ a Gaussian with the help
of a symbolic computation  software, one f\/inds complicated expressions
involving Bessel functions.

The distance between two arbitrary states on $\R^m$ is less easily computable. Let us f\/irst recall
what happens in the one dimensional case (cf.~e.g.~\cite{EG99}).

\begin{proposition}
For any two states $\varphi_1$ and $\varphi_2$,
whose corresponding probability measures $\mu_1$ and $\mu_2$
have no singular continuous part,
\[
\dc(\varphi_1,\varphi_2)=\int_{-\infty}^{+\infty}\left| \int_{-\infty}^{z}
\bigl(\de\mu_1(\xi) - \de\mu_2(\xi)\bigr)\right|\de z  .
\]
\end{proposition}

\begin{proof}
Integrating by parts one f\/inds
\[
\varphi_1(f)-\varphi_2(f)=-\int_{\R}f'(z)\Delta(z)  \de z,
\]
where
\begin{gather}
\Delta(z) \doteq \int_{-\infty}^{z}\big(\de\mu_1(\xi) - \de\mu_2(\xi)\big)
\label{eq:60}
\end{gather}
is the \emph{cumulative distribution} of the measure $\mu_1-\mu_2$.
We are assuming that $\mu_i$ are the sum of a~pure point part
and a part that is absolutely continuous with respect to the Lebesgue measure.
In this case $\Delta(z)$ is piecewise continuous, its sign is a piecewise continuous
function, and the primitive $h(z)$ of the sign of $\Delta(z)$ is a Lipschitz continuous
function.
Now, the $1$-Lipschitz condition says that a.e.~$|f'(z)|\leq 1$ on $\R$, hence
\[
\dc(\varphi_1,\varphi_2)\leq \int_{\R}|\Delta(z)| \de z.
\]
This upper bound is attained by any $1$-Lipschitz function $h$ such that
\[
h'(z) = 1 \quad \text{when} \ \ \Delta(z) \geq 0,\qquad h'(z) = -1 \quad \text{otherwise}.
\]
By the above consideration, such a $1$-Lipschitz function (at least one) always exists.
\end{proof}

In a quantum context, one can view the cumulative distributions $c_i = \int_\infty^z \de\mu_i(\xi)$
as the probability to f\/ind the particle on the half-line $(-\infty,z]$ before the transport ($i=1$) or after the transport ($i=2$). $\Delta(z)$ in~\eqref{eq:60} measures the probability f\/low across $z$, and the Wasserstein distance is the integral on $\R$ of the modulus of this probability f\/low.

On $\R^m$ with $m>1$ there is no
such an explicit result. If $\psi_1$,
$\psi_2$ are bounded Lipschitz functions with compact support, we
know from \cite{EG99} that there exist two (a.e.~unique) bounded measurable
Lipschitz functions $a,u:\R^m
\to\R$ such that $a\geq 0$,
\[
-\nabla(a\nabla u)=\psi_1-\psi_2
\]
in the weak sense, and $|\nabla u|=1$ almost everywhere on the set where $a>0$.
We have then
\[
\dc(\varphi_1,\varphi_2)=\int_{\R^m}a(x)\de x  .
\]
Indeed integrating by parts, we can write
\[
\varphi_1(f)-\varphi_2(f)=\int_{\R^m}(\nabla f)\cdot (a\nabla u)\de x
\]
and since $|a\nabla u|=a$ we have
\[
\dc(\varphi_1,\varphi_2)\leq \int_{\R^m}a(x)\de x \;.
\]
The sup is attained by the function $f=u$.

\section{Noncommutative examples}\label{sec:4}

At the light of Proposition~\ref{prop1} one may wonder if the analogy between the spectral and the Wasserstein distances still makes sense in a non-commutative framework.
In other terms, for~$\A$ noncommutative is the distance on $\Sc(\A)$ computed by~(\ref{eq:5.1}) related to some Wasserstein distance?

The most obvious answer, based on Gelfand's identif\/ication~(\ref{eq:16})  between points and pure states, would be to consider the Wasserstein distance  $W_D$ on the metric space $(\pa, d_D)$, and question whether $W_D$ on the set $\text{Prob}(\pa)$ of probabilities distributions
on $\pa$ coincides with $d_D$ on $\Sc(\A)$.
This is obviously true for pure states since by (\ref{eq:43})
\[
W_D(\omega_1,\omega_2) = c(\omega_1, \omega_2) =
d_D(\omega_1, \omega_2) \qquad \forall \; \omega_1, \omega_2 \in \pa.
\]
For non-pure states however this is usually not true. The reason is that even if
 \begin{gather*}
\Sc(\A)\subset \text{Prob}(\pa)
\end{gather*}
for commutative and almost commutative $C^*$-algebras (see the def\/inition in the next paragraph), there is not a 1-to-1 correspondence
between the two sets (except in the commutative case).
For instance, as recalled in Section~\ref{moyalsection}, $S(M_2(\C))$ is a non-trivial quotient of $\text{Prob}(\mathcal{P}(M_2(\C)))$, mea\-ning that two probability distributions $\phi_1\neq \phi_2$ may give the same state $\varphi_1=\varphi_2$. Thus, given a~spectral triple on $M_2(\C)$ as the one studied in~\cite{Moyal}, one has $d_D( \varphi_1, \varphi_2) = 0$ while $W_D(\phi_1, \phi_2) \neq 0$ since the cost being a distance
(and assuming $\pa$ is a polish space), $W$ is a distance and vanishes if and only if $\phi_1 =\phi_2$.

However this does not mean that the optimal transport interpretation of the spectral distance loses all interest in the noncommutative framework. As explained in Section~\ref{smsection}, there are interesting analogies between the spectral distance and some Wasserstein distances other than~$W_D$  in almost-commutative geometries. The latest are spectral triples $(\A, \HH, D)$ obtained as the product of a spin manifold by a f\/inite-dimensional spectral triple $(\A_I, \HH_I, D_I)$, namely
\begin{gather}
  \label{eq:55}
  \A  =\coi \otimes \A_I,\qquad \HH = L_2(\M, \Sc)\otimes \HH_I,\qquad D= -i\gamma^\mu\partial_\mu \otimes \III_I + \gamma^{m+1}\otimes D_I,
\end{gather}
where $\III_I$ is the identity of $\HH_I$ and $\gamma^{m+1}$ is the product of the gamma matrices in case $m$ is even, or the identity in case $m$ is odd. $\A_I$ being f\/inite-dimensional (or, equivalently, $\coi$ being Abelian) one has \cite{Kad86}:
\begin{gather}\label{eq:59}
   \mathcal{P}(\A) = \mathcal{P}(\coi) \times \mathcal{P}(\A_I)  .
\end{gather}
For simplicity we discuss here the case $\A_I = \C^2$ acting on $\HH_I= \C^2$, while in Section~\ref{smsection} we will focus on the f\/inite-dimensional spectral triple describing the internal degrees of freedom (hence the subscript $I$) of the standard model of particle physics.
Since $\C^2$ has two pure states~-- $\delta_0(z_0,z_1) =z_0$, $\delta_1(z_0,z_1) = z_1$ $\forall\; (z_0, z_1)\in \C^2$~-- from \eqref{eq:59}
\begin{gather*}
   \mathcal{P}(\A) = \M \times \{0,1\}
\end{gather*}
This is Connes' idea of ``product of the continuum by the discrete''  \cite{Con94} seen at the level of pure states: through the product by a f\/inite-dimensional spectral triple, the points of the manifold acquire a $\mathbb{Z}_2$ internal discrete structure, namely to a point $x=\delta_x$ in the commutative case corresponds two pure states in $\pa$,
\[
x_0 = (\delta_x,\delta_0),\qquad x_1 = (\delta_x,\delta_1) .
\]
Moreover, although $\pa$ is the disjoint union of two copies of $\M$, points on distinct copies need not to be at inf\/inite spectral distance from one another. Indeed  one shows that \cite{MW02}
\begin{gather}
d_D(x_0, x_1) = d_{D_I} (\delta_0, \delta_1),\label{eq:25}
\end{gather}
where $d_{D_I}$ denotes the spectral distance associated to $(\A_I, \HH_I, D_I)$. $\A_I$ is represented on $\HH_I$ by diagonal matrices but $D_I\in M_2(\C)$ needs not to be diagonal. Especially, if
$D_I$ has non-zero of\/f-diagonal terms then $\norm{[D,a]}\leq 1$ is a non-trivial constraint which guarantees that \eqref{eq:25} is f\/inite. Note that if one takes the direct sum of spectral triples instead of a product, namely $\coi \oplus \coi$ acting on $L_2(\M,\Sc)\oplus L_2(\M,\Sc)$, one gets the same pure states as above  but $D= -i\gamma^\mu\partial\mu \oplus  -i\gamma^\mu\partial\mu$ does not have of\/f-diagonal terms and the two copies of $\M$ are at inf\/inite distance from one another (points have not enough space to ``talk to each other'' through the of\/f-diagonal terms of $D$).

 From an optimal transport point of view, the non-vanishing of \eqref{eq:25} can be interpreted as the fact that ``staying at a point'', which is costless in the commutative case since $c(x,x) = d_D(x,x) = 0$, may have a cost
\begin{gather}
d_D(x_0,x_1)\neq 0\label{eq:61}
\end{gather}
 in almost-commutative geometries,
corresponding to the ``internal jump'' from $x_0$ to $x_1$. We investigate this idea in Section~\ref{smsection}, showing in Proposition~\ref{propsm} that the spectral distance between pure states is the minimal work $W_I$ associated to a cost~$c_I$.

\subsection{The Moyal plane}\label{moyalsection}

\begin{figure}[t]
\centerline{\includegraphics[scale=0.95]{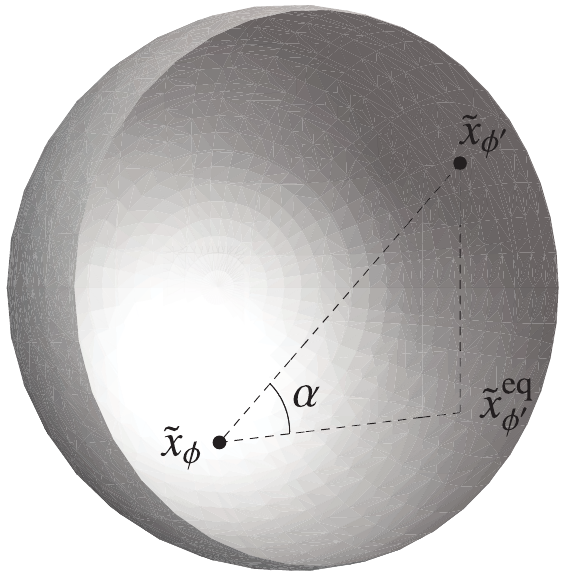}}

\caption{The vertical plane containing ${\bf \bar x}_{\phi}$, ${\bf \bar x}_{\phi'}$.}\label{Moyalfig}
\end{figure}

The Moyal algebra $\A_\theta$ is the noncommutative deformation of the non-unital Schwartz algebra $S(\R^2)$
with point-wise product
\[
(f\star g)(x)=\frac{1}{(\pi\theta)^2}\int d^2y d^2z\, f(x+y)g(x+z)e^{-i 2y \Theta^{-1}z} \qquad \forall \; f,g \in S(\R^2),
\]
where $y \Theta^{-1}z\equiv y^\mu \Theta^{-1}_{\mu\nu}z^\nu$ and
\[
 \Theta_{\mu\nu}=\theta\begin{pmatrix} 0&1 \\ -1& 0 \end{pmatrix}
\]
with $\theta\in\mathbb{R}$, $\theta\ne0$. In \cite{Moyal} we studied the spectral distance associated to the spectral triple built in \cite{Gayral:2004rc} around the action of $\A_\theta$ on $L_2(\R^2)$ and the usual Dirac operator on $\R^2$. We found that the topology induced by the spectral distance on $\mathcal{S}(\A_\theta)$ is not the weak*-topology (a condition required by Rief\/fel \cite{Rie99,Rie03,Rie04} in the unital case in order to def\/ine compact quantum metric spaces, and adapted by~\cite{Lat05} to the non-unital case). However by viewing $\A_\theta$ as an algebra of inf\/inite-dimensional matrices, we proposed some f\/inite-dimensional truncations of the Moyal spectral triple, based on the algebra
$M_n(\C)$, $n\in\N$, that makes it a quantum metric space. Explicitly, for $n=2$ $\mathcal{P}(M_2(\C))$ is homeomorphic to the Euclidean $2$-sphere,
\[
\xi \in \mathcal{P}(M_2(\C)) \,\longrightarrow\, {\bf x}_\xi = (x_\xi, y_\xi, z_\xi) \in S^2,
\]
so that a non-pure state $\varphi$ is determined by a probability distribution $\phi$ on $S^2$. Its evaluation
\[
    \varphi(a) = \int_{S^2} \phi({\bf x}_\xi )\; \text{Tr} (s_\xi  a) \,d {\bf x}_\xi  = \text{Tr}\left( \left(\int_{S^2} \phi({\bf x}_\xi ) s_\xi  \,d {\bf x}_\xi \right)  a\right) \qquad  \forall\; a\in M_2(\C),
\]
with $s_\xi\in M_2(\C)$ the support of the pure state $\mu^{-1}({\bf x}_\xi )$ and $d {\bf x}_\xi $ the $SU(2)$ invariant measure on $S^2$,  only depends on the barycenter of $\phi$,
\[
{\bf\bar x}_\phi = (\bar x_\phi,  \bar y_\phi,  \bar y_\phi) \qquad \text{with}\quad
 \bar x_\phi  := \int_{S^2} \phi({\bf x}_\xi) x_\xi d {\bf x}_\xi
\]
and similar notation for $\bar y_\phi$, $\bar z_\phi$. With the equivalence relation
$\phi\sim \phi' \Longleftrightarrow {\bf\bar x}_\phi = {\bf \bar x}_{\phi'}$
one gets that
\[
{\mathcal{S}}(\M_2(\C)) = {\mathcal{S}}(C(S^2))/\sim
\]
is homeomorphic to the Euclidean 2-ball:
\[
    \varphi \overset{\mu}{\longrightarrow} {\bf \bar x}_\phi \in {\mathcal B}^2 .
\]
In \cite{Moyal} we computed
\[
d_D({\bf \bar x}_\phi,{\bf \bar  x}_{\phi'})=
\sqrt{\frac{\theta}{2}}\times\begin{cases}
\cos\alpha\;  d_{Ec} (\tdx_\phi, \tdx_{\phi'}) & \mathrm{when }\ \alpha\leq \dfrac{\pi}4,\vspace{1mm}\\
\dfrac{1}{2\sin \alpha}\, d_{Ec} (\tdx_{\phi}, \tdx_{\phi'})  & \mathrm{when }\  \alpha\geq \dfrac{\pi}4,
\end{cases}
\]
where $d_{Ec} ({\bf \bar x}_{\phi}, {\bf \bar x}_{\phi'})=|{\bf \bar x}_{\phi}-{\bf \bar x}_{\phi'}|$ is the Euclidean distance
and $\alpha$ is the angle between the segment $[{\bf \bar x}_{\phi}, {\bf \bar x}_{\phi'}]$ and the horizontal plane $z_\xi = \text{const}$ (see  Fig.~\ref{Moyalfig}).

The state space of the truncated Moyal algebra is not the space of distributions on~$S^2$, but a quotient of it. Therefore  for $\phi\sim \phi', \phi\neq \phi'$, one has $d_D(\varphi, \varphi') = 0$ while $W_D(\varphi, \varphi') \neq 0$. In this example there seems to be no Wasserstein distance naturally associated to the spectral distance. In the next section, we exhibit another noncommutative example for which there exists a Wasserstein distance $W'$,  with cost $d'$ def\/ined on a set bigger than~$\pa$, that coincides with~$d_D$ for some non-pure states.

\subsection{The standard model}\label{smsection}

The spectral triple describing the standard model of elementary particles is the product \eqref{eq:55}  of a~manifold with
the f\/inite-dimensional algebra $\C \oplus \HHH \oplus M_3(\C)$ acting on $\HH_I = \C^{96}$.  $D_I$ is a $96\times 96$ matrix with entries the masses of the elementary fermions, the Cabibbo--Kobayashi--Maskawa matrix and the neutrino mixing-angles. The choice of the algebra is dictated by physics~\cite{Chamseddine:2007oz} (its unitaries group gives back the gauge group of the standard model) and 96 is the number of elementary fermions\footnote{(6 leptons $+$ 6 quarks $\times$ 3 colors) $\times$ 2 chiralities $=$ 48,
to which are added 48 antiparticles.}. The pure states of $M_3(\C)$ turn out to be at inf\/inite distance from one another and from any other pure state~\cite{MW02} so that, from the metric point of view, the interesting part is the product $(\A, \HH, D)$ of a manifold by ($\A_I \doteq \C \oplus \HHH, \HH_I, D_I)$. By \eqref{eq:59} the pure-states space is
\[
{\mathcal P} (\A) = \M \times\left\{0,1\right\}
\]
since both $\C$ and the algebra of quaternion $\HHH$ have one single pure state\footnote{$\HHH$ is a real, but not a complex algebra. Real $C^*$-algebras are def\/ined similarly to complex ones \cite{Goodearl:1982fk}, except that one
imposes that $1+a^*a$ is invertible~-- which comes as a consequence in the complex case. A state $\varphi$ is then def\/ined as
a real, real-linear, positive form such that $\varphi(\mathbb{I})=1$ and $\varphi(a^*) = \varphi(a).$ Hence $\HHH$ has only one state, given by half the trace. Viewing $\C$ and $\coi$ as real $C^*$-algebras, their states are obtained by taking the real part of their ``usual'' complex states. See \cite[page~42]{MarPhD}.}
\[
  \delta_{\C}(z) =\Re(z) \quad\forall\;  z\in\C,\qquad   \delta_{\HHH}(h) = \frac 12  \text{Tr} \,h \quad \forall\; h= \left(\begin{array}{cc} \alpha & \beta \\ -\bar\beta &\bar \alpha\end{array}\right) \in\HHH, \alpha, \beta\in\C.
\]
The space underlying the standard model is thus the disjoint union of two copies of the manifold, and with some computation one shows \cite{MW02} that the spectral distance is f\/inite on $\pa$ and coincides with the geodesic distance $d'$ on the manifold $\M' = \M\times \mathcal{I}$,
with $\mathcal{I}:=[0,1]$ the closed unit interval, with metric
\begin{gather}
  \label{eq:26}
  \left(\begin{array}{cc}
g^{\mu\nu}(x)& 0 \\ 0 &\norm{D_I}^2 \end{array}\right),
\end{gather}
where $g$ is the Riemannian metric on $\M$.
Assuming that $\M$ is complete, and writing $t\in\mathcal{I}$ the extra-coordinate, this can be restated as:
the spectral distance between pure states in the standard model coincides with the Wasserstein distance $W'$ on the metric space $(\M', d')$ restricted to the two hyperplanes $t=0$, $t=1$. This reformulation of the main result of \cite{MW02} allows to determine the spectral distance between a certain class of non-pure states, namely those localized on one of the copies of $\M$.

Explicitly, $\mathcal{S}(\A)$ is the set of couples of measures $\varphi =(\mu,\nu)$ on $\M$, normalized to
\[
  \int_\M   d\mu + \int_\M    d\nu = 1,
\]
whose evaluation on
\[
\A\ni a= f \oplus \left(\begin{array}{cc} g & b \\ -\bar b & \bar g\end{array}\right) ,
\]
where
$f$, $g$, $b$ are in $\coi$, is given by
\[
 \varphi (a) = \int_\M \Re(f) \,d\mu + \int_\M \Re(g) \,d\nu.
\]
Two states $\varphi_1$, $\varphi_2$ are localized on the same copy of $\M$ if $\varphi_1 = (0,\nu_1)$, $\varphi_2=(0,\nu_2)$; or
$\varphi_1=(\mu_1,0)$, $\varphi_2=(\mu_2,0)$. Such states can be viewed as elements of $\Sc(\A)$, $\Sc(C^\infty_0(\M'))$ and $\Sc(\coi)$.

\begin{proposition} \label{stateloc} For two states $\varphi_1$, $\varphi_2$ localized on the same copy of $\M$,
  \begin{gather}
d_D(\varphi_1, \varphi_2) = W(\varphi_1, \varphi_2) = W'(\varphi_1, \varphi_2).\label{eq:20}
\end{gather}
  \end{proposition}

\begin{proof}
To f\/ix notation we assume that $\varphi_1 = (0,\nu_1)$, $\varphi_2 = (0,\nu_2)$ are localized on the  $\delta_\HHH$ copy of $\M$, that we write $\M_\HHH$ and associate to the value $t=1$ of the extra-parameter $t\in \mathcal{I}$. The evaluation of $\varphi_1 -\varphi_2$ on $a\in\A$ only depends on $\text{Re}(g)$. Moreover for $a=a*$, $\norm{[D,a]}\leq 1$ implies that $g$ is $1$-Lipschitz \cite[equation~(16)]{MW02}, so that
\begin{gather}
d_D(\varphi_1, \varphi_2 )\leq W(\varphi_1, \varphi_2 ).\label{eq:63}
\end{gather}
The equality is attained by considering $a=g_w \otimes \mathbb{I}_3$  where $g_w$ is the real $1$-Lipschitz function that attains the supremum in the computation of $W$. Hence the f\/irst equation in~\eqref{eq:20}.

To show that $W(\varphi_1, \varphi_2 )= W'(\varphi_1, \varphi_2 )$, one f\/irst notices that~\eqref{eq:26} being block-diagonal implies that $d'((x, 1), (y, 1)) = d(x, y)$. Hence  the restriction $g(y) \doteq g'(y,1)$ on $\M_\HHH$ of any $1$-Lipschitz function $g'$ in $C^\infty_0(\M')$ is $1$-Lipschitz,
\[
 \abso{g'(x, t_0) - g'(y, t_1)}\leq d'((x,t_0), (y,t_1)) \Rightarrow  \abso{g(x) - g(y)}\leq d(x,y).
\]
Conversely to any $1$-Lipschitz function $g\in\coi$ one associates
\[
g'(x,t) \doteq g(x) \qquad\forall\; x\in\M, t\in \mathcal{I}\label{eq:65}
\]
which is $1$-Lipschitz in $C^\infty_0(\M')$ and takes the same value as $g$ on $\M_\HHH$. Therefore
\[
   \underset{g'\in C^\infty_0(\M'), \norm{g'}_{\text{Lip} =1}}{\sup} \abso{\varphi_1(g') -\varphi_2(g')} =    \underset{g\in C^\infty_0(\M), \norm{g}_{\text{Lip} =1}}{\sup} \abso{\varphi_1(g) -\varphi_2(g)}
\]
which, together with \eqref{eq:63},  yields the result.
\end{proof}

It is tempting to postulate that $d_D(\varphi_1, \varphi_2) = W'(\varphi_1, \varphi_2)$ for states localized on dif\/ferent copies, or states that are not localized on any copy (i.e.\  $\varphi =(\mu,\nu)$ with $\mu$, $\nu$ both non-zero). This point is under investigation.

A disturbing point in the computation of the spectral distance in the standard model is the appearance of a compact extra-dimension $\mathcal{I}$ while the internal structure is discrete.  In~\cite{MW02} $\mathcal{I}$  came out more
as a computational artifact than a requirement of the model. From the Wasserstein distance point of view, the introduction of the extra-dimension can be seen in the following proper inclusions:
\[
\pa \subset \M'\subset \sa.
\]
Namely $W'$ is associated to the metric space $(\M',d')$ which is bigger than $\pa$ (the points between the sheets are not pure states of $\A$) and smaller than $\sa$ (non-pure states localized on a copy are not in $\M'$). To get rid of the extra-dimension, one could consider the metric space $(\pa, d_D)$ with associated Wasserstein distance $W_D$, but this would be of poor interest since the def\/inition of $d_D = d'$ requires the knowledge of~$\M'$. Alternatively one could look for a cost def\/ined solely on $\M$. For states (pure or not) localized on the same copy, this cost is simply the geodesic distance on $\M$, as shown in Proposition~\ref{stateloc}. For pure states on distinct copies such a~cost $c_I$ also exists and is given by
\begin{gather}    \label{eq:41}
c_{I}(x, y) \doteq \sqrt{d(x,y)^2 + \frac 1{\norm{D_I}^2}}.
\end{gather}
Note that $c_I$ is not a distance since it does not vanish on the diagonal,
\[
  c_I(x,x) = \frac 1{\norm{D_I}}\qquad \forall \; x\in\M,
\]
but gives back the jump-cost (\ref{eq:61}) (by (\ref{eq:26}), $d(x_0,x_1) = \frac 1{\norm{D_I}}$). However $c_I$ satisf\/ies the condition required to proved Kantorovich's duality (namely \cite{Vil03} it is lower semicontinuous and satisf\/ies $c_I(x,y)\geq a(x) + b(y)$ for some real-valued upper semicontinous $\mu_i$-integrable functions).
Hence it makes sense to
consider the minimal work $W_I$ associated to $c_I$, whose formula is given by~(\ref{eq:8}) with the Lipschitz condition replaced by~$\abso{f(x) - f(y)} \leq c_I(x,y)$.
\begin{proposition}\label{propsm}
The spectral distance between pure states
 $x_0 \doteq (\delta_x, \delta_\C)$, $y_1\doteq (\delta_y, \delta_\HHH)$ on distinct copies
 is
\[
d_D(x_0, y_1) = W_I(x,y) .
\]
\end{proposition}

\begin{proof}
For pure states the minimal work is the cost itself: $W_I(x,y) = c_I(x,y)$ for any $x,y\in\M$. Since $d_D(x_0, y_1) = d'(x_0,y_1)$ as recalled above \eqref{eq:26}, the result follows if one proves that
  \begin{gather}
    \label{eq:39}
     d'^2(x_0,y_1) =  d^2(x,y)+ \frac 1{\norm{D_I}^2}.
  \end{gather}
This has been shown in \cite{MW02} but we brief\/ly restate the argument here for sake of completeness. Let us write $x^a= (x^\mu\in\M, t\in [0,1])$ a point of $\M'$ and $g^{tt} = g_{tt}^{-1} = \norm{D_I}^2$ the extra-metric component. The Christof\/fel symbols involving $t$ are
\[
\Gamma^t_{t\mu}=\Gamma^t_{\mu t}= \frac{1}{2}g^{t}\partial_{\mu} g_{t} ,\qquad \Gamma^\mu_{t}=-\frac{1}{2} g^{\mu\nu} \partial_\nu g_{tt} ,
\qquad \Gamma^{\mu}_{0\nu}=\Gamma^{\mu}_{\nu t} = \Gamma^{t}_{tt}=\Gamma^t_{\mu\nu}=0.
\]
The geodesic equation $\ddot x^a + \Gamma_{bc}^a \dot x^b \dot x^c$ writes
\begin{subequations}
\begin{gather}
\ddot x^t + g^{tt}(\partial_\mu g_{tt})
\dot x^t \dot x^{\mu}= 0 ,
\\ \label{geozero}
\ddot x^\mu - \frac{1}{2} g^{\mu\nu}(\partial_\nu g_{tt})
(\dot t)^2 +\Gamma^{\mu}_{\lambda\rho}\dot x^\lambda \dot x^\rho = 0
\end{gather}
\end{subequations}
and, because $g_{tt}$ is a constant, they simplify to
\begin{subequations}
\begin{gather}
\dot t= \text{const} = K  , \\
\ddot x^\mu+\Gamma^{\mu}_{\lambda\rho}\dot x^\lambda \dot x^\rho = 0.\label{geo4}
\end{gather}
\end{subequations}
The f\/irst geodesic equation indicates that the  proper length $\tau$ of a geodesic $\mathcal{G}' = x(\tau)$ in $\M'$ between $x_0$ and $y_1$ is proportional to $t$,
\[
  dt= K d\tau,
\]
as well as to the line element $ds = \sqrt{g_{\nu\nu}dx^\mu dx^\nu}$ of $\M$
since
\begin{gather}\label{dessin}
1= \norm{\dot x} = g_{ab} \dot x^a \dot x^b =
g_{\mu\nu}\frac{dx^\mu}{d\tau}\frac{dx^\nu}{d\tau}+ g_{tt}K^2 = \frac{ds^2}{d\tau^2} + g_{tt}K^2,
\end{gather}
so that, assuming $g_{tt}K^2\neq1$ (which from \eqref{dessin} amounts to take $x\neq y$),
\[
d\tau=  \frac{ds}{\sqrt{1 - g_{tt}K^2}}\;.
\]
The second geodesic equation~\eqref{geo4} shows that  the projection on $\M$ of
$\mathcal{G}'$ is a geodesic~$\mathcal{G}$ of~$\M$. Therefore
\[
  d'(x_0, y_1) = \int_{\mathcal{G}'}  d\tau = \int_{\mathcal{G}} \frac{d\tau}{ds} ds =  \frac 1{\sqrt{1 - g_{tt}K^2}}\int_{\mathcal{G}} d s=  \frac 1{\sqrt{1 - g_{tt}K^2}} d(x,y),
\]
that is to say
\[
  d'^2(x_0, y_1) = d^2(x,y) + g_{tt}K^2   d'^2(x_0, y_1).
\]
Writing $K^2 d'^2(x_0,y_1)$ as
\[
\left( \int_{\mathcal{G'}} K d\tau\right)^2 = \left(\int_{\mathcal{G'}} dt\right)^2 = \left(t((y,1) - t((x,0))\right)^2 = 1
\]
one obtains~\eqref{eq:39}, and the result.
\end{proof}

 Let us underline an interesting feature of this proposition, namely that a cost which is not a~distance $\M$ can be seen as a distance on $\M\times \M$.

Proposition \ref{propsm} can also be rewritten using  a cost vanishing on the diagonal, namely
\[
  c'_I(x,y) \doteq c_I(x,y) - \frac 1{\norm{D_I}^2} = \sqrt{d(x,y)^2 + \frac 1{\norm{D_I}^2}} - \frac 1{\norm{D_I}} .
\]
Writing $W_I'$ the associated minimal work, one has
\[
  d_D(x_0, y_1) = W_I(x,y) + \frac 1{\norm{D_I}^2}.
\]
Quite remarkably, the cost \eqref{eq:41} is similar to the cost~(32) introduced in~\cite{Brenier} in the framework of the relativistic heat equation.

Proposition~\ref{propsm} relies on the fact that the jump-cost~\eqref{eq:61} is constant. From a physics point of view, this means that one does not take into account the Higgs f\/ield.
In almost commutative geometries, the latest is obtained by \emph{inner fluctuation of the metric} \cite{Con96} that substitute $D$ with a covariant Dirac operator. From the metric point of view this
amounts to replacing (see \cite{MW02} for details)
 \eqref{eq:26} by
\[
  \left(\begin{array}{cc}
g^{\mu\nu}(x)& 0 \\ 0 &\norm{D_I +H(x)}^2 \end{array}\right),
\]
where
\[
H(x) = \left(\abso{1+h_1(x)}^2 + \abso{h_2(x)}^2\right)m_t^2
 \]
 where $(h_1,h_2)$ is the (complex) Higgs doublet and~$m_t$ is the mass of the quark top. Instead of~\eqref{eq:25} one has $d_D(x_0, x_1) = d_{D_I + H(x)}(\delta_\C, \delta_\HHH )$ (the jump-cost is no longer constant). In analogy with~(\ref{eq:41}), one could def\/ine the cost
\[
  \tilde c_I(x,y) = \sqrt{d(x,y)^2 + d_{D_I + H(x)}(\delta_\C, \delta_\HHH )^2}
\]
which allows to avoid the introduction of the extra-dimension. However the projection of a~geodesic in~$\M'$ is no longer a geodesic of~$\M$ (\eqref{geozero} does not simplify to~\eqref{geo4}) so there is no way to express $d'(x_0, y_1)$
as a function of $d(x,y)$, and $d_D(x_0, y_1)= d'(x_0, y_1)$ no longer equals $W_I(x_0,y_1)$.

\section{Conclusion}
The spectral distance between states on a complete Riemannian spin mani\-fold coincides with the Wasserstein distance of order $1$. In the noncommutative case  the analogies between the spectral distance and various Wasserstein distances, although still not fully understood, shed an interesting light on the interpretation of the distance formula in noncommutative geometry. In physics, def\/ining the distance as a supremum rather than an inf\/imum is useful since, at small scale, quantum mechanics indicates that notions as ``paths between points'' no longer make sense, so that the classical def\/inition of distance as the length of the shortest path loses any operational meaning. An interesting feature of noncommutative geometry is to provide a notion of distance overcoming this dif\/f\/iculty. In transport theory the interpretation of Kantorovich formula has an interpretation in economics rather than in quantum mechanics: while Monge formulation corresponds to the minimization of a cost, Kantorovich dual formula corresponds to the maximization of a prof\/it. To repeat a classical example found in the literature:
assume that the distribution of f\/lour-producers on a given territory~$\M$ is given by~$\mu_1$ and the distribution of bakeries by $\mu_2$. Consider a transport-consortium whose job consists in buying the f\/lour at factories and selling it to bakers. The consortium f\/ixes the value $f(x)$ of the f\/lour at the point~$x$ (it buys the f\/lour at the price $f(x)$ if there
is a factory, or it sells it at a~price $f(x)$ if there is a~bakery). The Wasserstein distance is the maximum prof\/it the consortium may hope, under the constraint of staying competitive, that means  not selling the f\/lour to a~price higher than the bakeries would pay if they were doing the transport by themselves (i.e.~$f(x)\leq f(y) + c(x,y)$ for all~$x$,~$y$).
This raises an interesting question: in a quantum context, what is the physical meaning of this ``prof\/it'' that one is maximizing while computing the distance?  If one view the states on~$\M$ as wave functions,
is the distance related to the minimum work required to transform one wave into the second?
More specif\/ically,  for $\M=\R^m$ what physical quantity represents~(\ref{eq:r2dist}), and what is the meaning of the point~$\alpha$?

\subsection*{Acknowledgments}
We would like to thank Hanfeng Li and anonymous referees for their valuable comments. 

\pdfbookmark[1]{References}{ref}
\LastPageEnding


\begin{thebibliography}{99}

\footnotesize\itemsep=0pt

\bibitem{Amb00}
Ambrosio L.,
Lecture notes on optimal transport problems, in Mathematical Aspects of Evolving Interfaces (Funchal, 2000), {\it Lecture Notes in Math.}, Vol.~1812, Springer, Berlin, 2003, 1--52.

\bibitem{AFLR07}
 Azagra D., Ferrera J., L\'opez-Mesas F., Rangel Y.,
 Smooth   approximation of Lipschitz functions on Riemannian manifolds,
\href{http://dx.doi.org/10.1016/j.jmaa.2006.03.088}{{\it J. Math. Anal. Appl.}} \textbf{326} (2007), 1370--1378,
\href{http://www.arxiv.org/abs/math.DG/0602051}{math.DG/0602051}.

\bibitem{BR}
Bratteli O., Robinson D.W.,
Operator algebras and quantum statistical mechanics. 1.~$C^*$- and $W^*$-algebras, symmetry groups, decomposition of states,
2nd ed., {\it Texts and Monographs in Physics}, Springer-Verlag, New York, 1987.

\bibitem{BV01}
Biane P., Voiculescu D.,
A free probability analogue of the Wasserstein metric on the trace-state space,
\href{http://dx.doi.org/10.1007/s00039-001-8226-4}{{\it Geom. Funct. Anal.}} \textbf{11} (2001), 1125--1138,
\href{http://www.arxiv.org/abs/math.OA/0006044}{math.OA/0006044}.

\bibitem{bls}
Bimonte G., Lizzi F., Sparano G.,
Distances on a lattice from non-commutative geometry,
\href{http://dx.doi.org/10.1016/0370-2693(94)90302-6}{{\it Phys. Lett. B}} {\bf 341} (1994), 139--146,
\href{http://www.arxiv.org/abs/hep-lat/9404007}{hep-lat/9404007}.

\bibitem{Brenier}
Brenier Y.,
Extended Monge--Kantorovich theory, in Optimal Transportation and Applications (Martina Franca, 2001),
{\it Lecture Notes in Math.}, Vol.~1813, Springer, Berlin, 2003, 91--121.


\bibitem{Moyal}
Cagnache E., D'Andrea F., Martinetti P., Wallet J.C.,
The Spectral distance in the Moyal plane,
\href{http://www.arxiv.org/abs/0912.0906}{arXiv:0912.0906}.

\bibitem{Chamseddine:2007oz}
 Chamseddine A.H., Connes A., Marcolli M.,
Gravity and the standard model with neutrino mixing,
{\it Adv. Theor. Math. Phys.} {\bf 11} (2007), 991--1089,
\href{http://www.arxiv.org/abs/hep-th/0610241}{hep-th/0610241}.

\bibitem{Cha93}
Chavel I.,
Riemannian geometry~-- a modern introduction,
{\it Cambridge Tracts in Mathematics}, Vol.~108, Cambridge University Press, Cambridge, 1993.

\bibitem{Con89}
Connes A.,
Compact metric spaces, Fredholm modules, and hyperf\/initeness,
\href{http://dx.doi.org/10.1017/S0143385700004934}{{\it Ergodic Theory Dynam. Systems}} \textbf{9} (1989), 207--220.

\bibitem{Con94}
Connes A.,
Noncommutative geometry, Academic Press, Inc., San Diego, CA, 1994.

\bibitem{Con95}
Connes A.,
Noncommutative geometry and reality,
\href{http://dx.doi.org/10.1063/1.531241}{{\it J. Math. Phys.}} \textbf{36} (1995), 6194--6231.

\bibitem{Con96}
Connes A.,
Gravity coupled with matter and the foundation of non-commutative geometry,
\href{http://dx.doi.org/10.1007/BF02506388}{{\it Comm. Math. Phys.}}  \textbf{182} (1996), 155--176.

\bibitem{connesckm}
Connes A.,
A unitary invariant in Riemannian geometry,
\href{http://dx.doi.org/10.1142/S0219887808003284}{{\it Int. J. Geom. Methods Mod. Phys.}} {\bf 5} (2008), 1215--1242,
\href{http://www.arxiv.org/abs/0810.2091}{arXiv:0810.2091}.

\bibitem{CL92}
Connes A., Lott J.,
The metric aspect of noncommutative geometry,
in New symmetry principles in quantum f\/ield theory (Carg\`ese, 1991), {\it NATO Adv. Sci. Inst. Ser. B Phys.}, Vol.~295, Plenum, New York, 1992, 53--93.

 \bibitem{Connes:2008kx}
Connes A., Marcolli M.,
Noncommutative geometry, quantum f\/ields and motives, {\it American Mathematical Society Colloquium Publications}, Vol.~55, American Mathematical Society, Providence, RI; Hindustan Book Agency, New Delhi, 2008.

\bibitem{DFR94}
Doplicher S., Fredenhagen K., Roberts J.E.,
Spacetime quantization induced by classical gravity,
\href{http://dx.doi.org/10.1016/0370-2693(94)90940-7}{{\it Phys. Lett.~B}} \textbf{331} (1994), 39--44.

\bibitem{DFR95}
Doplicher S., Fredenhagen K., Roberts J.E.,
The quantum structure of spacetime at the Planck scale and quantum f\/ields,
\href{http://dx.doi.org/10.1007/BF02104515}{{\it Comm. Math. Phys.}} \textbf{172} (1995),  187--220,
\href{http://www.arxiv.org/abs/hep-th/0303037}{hep-th/0303037}.

\bibitem{EG99}
 Evans L.C., Gangbo W.,
 Dif\/ferential equations methods for the Monge--Kantorevich mass transfer problem,
{\it Mem. Amer. Math. Soc.} {\bf 137} (1999), no.~653.


\bibitem{Gayral:2004rc}
Gayral V., Gracia-Bond{\'\i}a J.M., Iochum B., Sch{\"u}cker T., Varilly J.C.,
Moyal planes are spectral triples,
\href{http://dx.doi.org/10.1007/s00220-004-1057-z}{{\it Comm. Math. Phys.}} {\bf 246}  (2004), 569--623,
\href{http://www.arxiv.org/abs/hep-th/0307241}{hep-th/0307241}.

\bibitem{GS84}
 Givens C.R., Shortt R.M.,
A class of Wasserstein metrics for probability distributions,
\href{http://dx.doi.org/10.1307/mmj/1029003026}{{\it Michigan Math. J.}} \textbf{31} (1984), 231--240.

\bibitem{Goodearl:1982fk}
Goodearl K.R.,
Notes on real and complex $C^*$-algebras, {\it Shiva Mathematics Series}, Vol.~5,
Shiva Publishing Ltd., Nantwich, 1982.


\bibitem{fgbv}
Gracia-Bond{\'\i}a J.M., Varilly J.C., Figueroa H.,
Elements of noncommutative geometry,
{\it Birkh\"auser Advanced Texts: Basler Lehrb\"ucher}, Birkh\"auser Boston, Inc., Boston, MA, 2001.

\bibitem{GW79}
Greene R.E., Wu H.,
$C^\infty$ approximations of convex, subharmonic, and plurisubharmonic functions,
\href{http://www.numdam.org/item?id=ASENS_1979_4_12_1_47_0}{{\it Ann. Sci. {\'E}cole Norm. Sup. (4)}} \textbf{12} (1979),  47--84.

\bibitem{Gromov:1999fk}
Gromov M.,
Metric structures for Riemannian and non-Riemannian spaces,
{\it Progress in Mathematics}, Vol.~152, Birkh\"auser Boston, Inc., Boston, MA, 1999.

\bibitem{IKM01}
Iochum B., Krajewski T., Martinetti P.,
Distances in f\/inite spaces  from noncommutative geometry,
\href{http://dx.doi.org/10.1016/S0393-0440(00)00044-9}{{\it J. Geom. Phys.}} \textbf{37} (2001), 100--125,
\href{http://www.arxiv.org/abs/hep-th/9912217}{hep-th/9912217}.

\bibitem{Kad86}
Kadison R.V., Ringrose J.R.,
Fundamentals of the theory of operator algebras. Vol.~II.~Advanced theory,
{\it  Pure and Applied Mathematics}, Vol.~100, Academic Press, Inc., Orlando, FL, 1986.

\bibitem{Kan42}
Kantorovich L.V.,
On the transfer of masses,
{\it Dokl. Akad. Nauk. SSSR} \textbf{37} (1942), 227--229.

\bibitem{KR58}
 Kantorovich L.V., Rubinstein G.S.,
 On a space of totally additive functions,
 {\it Vestnik Leningrad. Univ.} \textbf{13} (1958), no.~7, 52--59.

\bibitem{Lat05}
Latr{\'e}moli{\`e}re F.,
Bounded-Lipschitz distances on the state space of a $C^*$-algebra,
{\it Taiwanese J. Math.} \textbf{11} (2007),  447--469.

\bibitem{MarPhD}
 Martinetti P.,
Distances en g\'eom\'etrie non-commutative, PhD Thesis,
\href{http://www.arxiv.org/abs/math-ph/0112038}{math-ph/0112038}.

\bibitem{Mar06}
Martinetti P.,
Carnot--Carath\'eodory metric and gauge f\/luctuations in noncommutative geometry,
\href{http://dx.doi.org/10.1007/s00220-006-0001-9}{{\it Comm. Math. Phys.}} {\bf 265} (2006), 585--616,
\href{http://www.arxiv.org/abs/hep-th/0506147}{hep-th/0506147}.

\bibitem{Mar08}
Martinetti P.,
Spectral distance on the circle,
\href{http://dx.doi.org/10.1016/j.jfa.2008.07.018}{{\it J. Funct. Anal.}} \textbf{255} (2008),   1575--1612,
\href{http://www.arxiv.org/abs/math.OA/0703586}{math.OA/0703586}.

\bibitem{Mardev}
 Martinetti P.,
 Smoother than a circle or  how noncommutative geometry provides the torus with an ego\-cent\-red metric,
in Modern Trends in Geometry and Topology (Deva, 2005), Cluj Univ. Press, Cluj-Napoca, 2006, 283--293,
\href{http://www.arxiv.org/abs/hep-th/0603051}{hep-th/0603051}.


\bibitem{MW02}
Martinetti P., Wulkenhaar R.,
Discrete Kaluza--Klein from scalar f\/luctuations in noncommutative geometry,
\href{http://dx.doi.org/10.1063/1.1418012}{{\it J. Math. Phys.}} \textbf{43} (2002), 182--204,
\href{http://www.arxiv.org/abs/hep-th/0104108}{hep-th/0104108}.

\bibitem{Mon81}
Monge G.,
M\'emoire sur la Th\'eorie des D\'eblais et des Remblais,
  Histoire de l'Acad. des Sciences de Paris, 1781.

\bibitem{Nas54}
Nash J.,
$C^1$-isometric imbeddings,
\href{http://dx.doi.org/10.2307/1969840}{{\it Ann. of Math. (2)}} \textbf{60} (1954), 383--396.

\bibitem{Nas56}
Nash J.,
The imbedding problem for Riemannian manifolds,
\href{http://dx.doi.org/10.2307/1969989}{{\it Ann. of Math. (2)}} \textbf{63} (1956), 20--63.

\bibitem{Rie99}
 Rief\/fel M.A.,
 Metric on state spaces,
 {\it Doc. Math.} \textbf{4}  (1999), 559--600,
 \href{http://www.arxiv.org/abs/math.OA/9906151}{math.OA/9906151}.

\bibitem{Rie03}
Rief\/fel M.A.,
Compact quantum metric spaces,
in Operator Algebras, Quantization, and Noncommutative Geometry,
{\it Contemp. Math.}, Vol.~365, Amer. Math. Soc., Providence, RI, 2004, 315--330,
\href{http://www.arxiv.org/abs/math.OA/0308207}{math.OA/0308207}.

\bibitem{Rie04}
Rief\/fel M.A.,
Gromov--Hausdorf\/f distance for quantum metric spaces,
{\it Mem. Amer. Math. Soc.} \textbf{168} (2004), no.~796, 1--65,
\href{http://www.arxiv.org/abs/math.OA/0011063}{math.OA/0011063}.

\bibitem{Roe96}
Roe J.,
Index theory, coarse geometry, and topology of manifolds,
{\it CBMS Regional Conference Series in Mathematics}, Vol.~90, American Mathematical Society, Providence, RI, 1996.

\bibitem{Vil03}
Villani C.,
Topics in optimal transportation,
{\it Graduate Studies in Mathematics}, Vol.~58, American Mathematical Society, Providence, RI, 2003.

\bibitem{Vil08}
Villani C.,
Optimal transport. Old and new,
{\it Grundlehren der Mathematischen Wissenschaften}, Vol.~338, Springer-Verlag, Berlin, 2009.

\end{thebibliography}
\end{document}